% ...

\documentclass[]{elsarticlepreprint}

\usepackage[latin1]{inputenc}
\usepackage{amsmath}
\usepackage{amssymb}
\usepackage{graphicx} 
\usepackage{graphics}
\usepackage{algorithm}
\usepackage{algorithmic}

\newtheorem{LEM}{Lemma}[section]
\newtheorem{PROP}{Proposition}[section]

\input comment.sty
\excludecomment{estensione}  

\begin{document}

\begin{frontmatter}

\title{Parabolic inverse convection-diffusion-reaction problem solved using an adaptive parametrization}

\author[Padova]{G. Deolmi} 
\author[Padova]{F. Marcuzzi} 
\address[Padova]{Dipartimento di Matematica Pura ed Applicata, Università degli Studi di Padova, Via Trieste 63, 35131 Padova, Italy \\ gdeolmi@math.unipd.it, marcuzzi@math.unipd.it}

\begin{abstract}
This paper investigates the solution of a parabolic inverse problem based upon the convection-diffusion-reaction equation, which can be used to estimate both water and air pollution. We will consider both known and unknown source location: while in the first case the problem is solved using a projected damped Gauss-Newton, in the second one it is ill-posed and an adaptive parametrization with time localization will be adopted to regularize it. To solve the optimization loop a model reduction technique (Proper Orthogonal Decomposition) is used.
\end{abstract}


\end{frontmatter}

\section{Introduction}

Inverse heat or mass convection problems, classically deal with the estimation of wall heat flux densities or intensities of source terms \cite{moutsoglou,park,girault2004,mcgrail}. As mentioned in \cite{girault2008}, inverse problems are usually mathematically ill-posed and regularization methods have been developed to ensure stable solutions. For an overview, see \cite{woodbury,isakov,kirsch} for example. Classical methods are penalization such as Tikhonov's regularization \cite{tikhonov}, or Bayesian methods using prior information \cite{aster}, iterative regularization \cite{alifanov} and regularization using singular value decomposition followed by truncation of the singular values spectrum \cite{shenefelt}.

In this paper we are interested in solving an inverse convection problem, whose direct model coincides with a parabolic convection diffusion reaction equation on a fixed domain. To deal with its ill-posedness we adopt a regularization algorithm based upon Truncated Singular Value Decomposition (TSVD) and diagonal scaling \cite{nocedal}; moreover an adaptive parametrization with time localization is formulated, to reduce the computational cost of the Gauss Newton algorithm and to obtain a better conditioned sensitivity matrix (cfr. e.g. figure \ref{cond_adattativa}).

Convection-diffusion-reaction equation can be used to model a variety of physical problems. For example in \cite{hamdi}, this equation is used to predict water quality in rivers, by measuring the quantity of organic matter contained. 
The importance of these pollution is estimated by the measures of the so-called BOD (Biologic Oxygen Demand) and COD (Chemical Oxygen Demand).  
In \cite{hamdi} the problem of identifying the location and the magnitude (intensity) of pollution point sources from the measurements of BOD on a part of the river is considered: the problem of source term identifycation is solved using an algorithm based on the minimization of a cost function of Kohn and Vogelius type. 
Also in \cite{rap} water pollution is considered: knowing the \textit{origin} of the source of contamination is probably the most important aspect when attempting to understand, and therefore to control, the pollution transport process. Thus, a challenging issue in environmental problems is the \textit{identification of sources} of pollution in waters. \cite{rap} deals with source identification problem, using Boundary Elemet Method (BEM). In \cite{girault2008}, the same problem of source estimation is considered to estimate the time-varying emission rates of pollutant sources in a ventilated enclosure, assuming that the velocity field is stationary: in fact in the frame of occupational risk prevention, the knowledge of both space and time distributions of contaminant concentration is a crucial issue to evaluate the workers exposure. Althought air pollution is considered, instead of water, the underlying model is still a convection-diffusion-reaction equation, with a different convective velocity field. In \cite{girault2008} source's location is supposed to be known. Possible applications of this study are concerned with cartography of pollutants in buildings, estimation of contaminant emission rates inside manufactures, leak detection, environment and process control through 'intelligent sensors' (controlled ventilation with closed-loop function of pollution threshold). A similar problem is considered in \cite{bracconier}. Finally in \cite{davoin}, a convection inverse problem is solved to determine an estimate of the source term as a funcition of the altitude and the temporal of iodine-131, caesium-134 and caesium-137 in the Chernobyl disaster.

In general, in inverse convection problems, either distributed control \cite{girault2008,rap}, or boundary control \cite{zabaras} or both \cite{lube} are considered. In the present paper we are interested in estimated location and intensity of pollution, and we assume to deal with boundary control, i.e. we suppose that the sources are located along domain's boundary. Thus, as in \cite{zabaras}, we deal with an inverse problem in which one is looking for the unknown conditions in part of the boundary, while overspecified boundary conditions are supplied in another part of the boundary (here the outflow region). As mentioned above, this type of problem can model both water and air pollution. 

As mentioned e.g. in \cite{gunzburger}, in \textit{inverse problems} or \textit{optimal control} or \textit{optimization settings}, one is faced with the need to do multiple state solves during an iterative process that determines the optimal solution. If one approximates the state in the reduced, $k$-dimensional space and if $k$ is small, then the cost of each iteration of the optimizer would be very small relative to that using full, high-fidelity state approximations. Thus \textit{Proper Orthogonal Decomposition (POD)} will be adopted in this paper as model reduction technique, to bring our study closer to a real time problem.

In section \ref{diretto} the direct problem is described.
In section \ref{inverso}, the inverse problem is formulated. Section \ref{known} deals with the problem of known source location, while in section \ref{pa} also source position is estimated.

\section{Description of the direct problem}
\label{diretto}
Let $[0,t_f)\subset\mathbb{R}$ and $\Omega$ be an open, limited and Lipschitz continuous boundary subset $\Omega\subset\mathbb{R}^2$, sufficiently regular. We denote with $\partial\Omega$ the boundary of $\Omega$. Let $C:\ [0,t_f)\times\Omega \rightarrow \mathbb{R}$, $C=C(t,\textbf{x})$ be the solution of the following (direct) parabolic convection-diffusion-reaction equation: 
\begin{equation}
\left\{\begin{array}{r l l l}
\frac{\partial C}{\partial t}-\mu \Delta C + \nabla \cdot (\textbf{u} C)+\sigma C & = & 0, & \qquad in \quad (0,t_f)\times\Omega \\ [5 pt]
C & = & C_0, & \qquad on \quad \left\{0\right\}\times\Omega \\ 
C & = & C_{in}, & \qquad on \quad (0,t_f)\times\Gamma_{in} \\ 
C & = & C_{up}, & \qquad on \quad (0,t_f)\times\Gamma_{up} \\ 
\mu\frac{\partial C}{\partial n} & = & 0, & \qquad on \quad (0,t_f)\times\Gamma_{down}\\ 
C & = & 0, & \qquad on \quad (0,t_f)\times\Gamma_{r}
\end{array}\right.
\label{direct_problem}
\end{equation}
where $\Gamma_{in}$, $\Gamma_{up}$, $\Gamma_{down}$ and $\Gamma_{r}$ are given disjoint sets such that $\partial \Omega = \Gamma_{in} \cup \Gamma_{up} \cup \Gamma_{down} \cup \Gamma_{r}$. 

Suppose that $C_{in}\in H^{\frac{1}{2}}(\Gamma_{in})$, $C_{up}\in H^{\frac{1}{2}}(\Gamma_{up})$, the initial condition $C_0\in L^2(\Omega)$ and the coefficients are independent on time, moreover $\mu\in L^{\infty}(\Omega)$, $\mu(\textbf{x})\geq \mu_0>0$ for all $\textbf{x}\in\Omega$, $\sigma\in L^{\infty}(\Omega)$, $\sigma(x)\geq0\ \text{a.e. in }\Omega$, $\textbf{u}\in[L^{\infty}(\Omega)]^2$, $div (\textbf{u})\in L^2(\Omega)$ are known. The \textit{direct problem} consists in finding the concentration $C$ over $\Omega$ at time $t_f$. 

As in \cite{girault2008}, we assume that the physical properties of the fluid are constant and that the transported contaminant is considered as a passive scalar, which means that it does not affect the velocity field. Thus we suppose to know $\textbf{u}$. 

\begin{figure}[h]
\begin{center}
\includegraphics*[width=12cm]{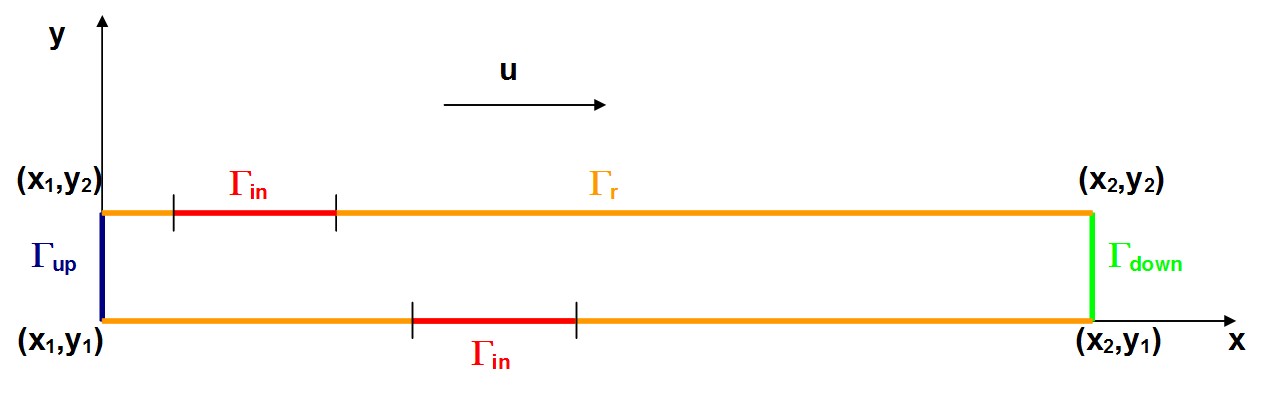}
\end{center}
\caption{\small \textit{Example of problem's domain $\Omega$.}}\normalsize
\label{dominio_fiume}
\end{figure}

An example of the 2D domain $\Omega$ is illustrated in figure \ref{dominio_fiume}.

\subsection{Wellposedness of the direct problem and finite element discretization}
Let $H^1_{\Gamma_r\cup\Gamma_{up}\cup \Gamma_{in}}(\Omega)$ be the set of $v\in H^1(\Omega)$ s.t. $v_{\left.\right|_{\Gamma_r\cup\Gamma_{up}\cup \Gamma_{in}}}=0$. Given $V\subset H^1_{\Gamma_r\cup\Gamma_{up}\cup \Gamma_{in}}(\Omega)$, the weak formulation of (\ref{direct_problem}) consists in finding $C\in L^2(0,t_f;H^1(\Omega))\cap\mathcal{C}^0 ([0,t_f);L^2(\Omega))$ s.t. 
\begin{equation}
\label{weak}
\begin{array}{r l l l}
\frac{d}{d t} (C(t),v) + a(\textbf{u}(t); C(t),v) & = & 0,& \forall v\in V,\\
C(0) & = & C_0,& in\ \Omega,
\end{array}
\end{equation}
where $a(\textbf{u}; \cdot,\cdot)$ is a bilinear form defined as 
$$a(\textbf{u}; w,v):= \int_{\Omega}k\nabla{w}\nabla v d\omega+ \int_{\Omega}\textbf{u} \cdot\nabla{w}v d\omega+\int_{\Omega}\sigma w v d\omega.$$

The wellposedness of the variational formulation is studied e.g. in \cite{quarteroni:Napde}. 

Consider now two families of subspaces $\left\{W_h,\ h>0\right\}$ and $\left\{V_h,\ h>0\right\}$ of $H^1(\Omega)$ and $V$ respectively, and let $C_{0,h}\in W_h$ be a suitable approximation of $C_0$. Then the Finite Element (FE) discretization of (\ref{weak}) consists in finding $C_h\in W_h$ s.t. 
\begin{equation}
\label{weak}
\begin{array}{r l l l}
\frac{d}{d t} (C_h(t),v_h) + a(\textbf{u}(t); C_h(t),v_h) & = & 0,& \forall v_h\in V_h,\\
C_h(0) & = & C_{0,h},& in\ \Omega.
\end{array}
\end{equation}

Given a basis of $W_h$, $\left\{\phi_i\right\}$, $i=1,\ldots,N_h$, where $N_h$ denotes the number of nodes in $\Omega$, it is well known that the FE discretization is equivalent to the solution of the following system of ODE's:
\begin{equation}
\begin{array}{r l l}
M \dot{\textbf{C}}(t)+ A(\textbf{u}(t)) \textbf{C}(t) & = & \textbf{F}(C_{in}),\\
\textbf{C}(0) & = & \textbf{C}_0.
\label{unreduced}
\end{array}
\end{equation}
where $M_{ij}=(\phi_i,\phi_j)$, $A(\textbf{u})_{ij}=a(\textbf{u}; \phi_i,\phi_j)$ and $\textbf{F}(C_{in})$ involves boundary conditions, in particular $C_{in}$. 

Given a time step $\Delta t$, consider a uniform subdivision of $[0,t_f)$ s.t. $(N-1) \Delta t = t_f$. Discretizing (\ref{unreduced}) in time, using e.g. the implicit euler method, we obtain
\begin{equation}
\begin{array}{r l l}
\left(M +\Delta t A(\textbf{u}(k+1))\right) \textbf{C}(k+1) & = & M \textbf{C}(k)+\Delta t \textbf{F}(C_{in}),\\
\textbf{C}(0) & = & \textbf{C}_0.
\end{array}
\end{equation}

\subsection{Proper Orthogonal Decomposition (POD) reduction}
\label{sezione_pod}
To obtain a faster solution algorithm, we adopt a reduction technique. A complete overview of all classical methods can be found e.g. in \cite{antoulas, schilders}. Since the right hand side in (\ref{unreduced}) depends on boundary conditions, it varies at each iteration of the optimization loop used to solve the inverse problem. As a consequence techniques largely used for linear constant matrices problems, like e.g. Balanced Truncation (BT), becomes too costly to be used. Thus we choose to adopt the Proper Orthogonal Decomposition (POD) method: althought its basis is stricly related to local dynamics, it is less costly to compute. In this paper we are focusing on the inverse problem solution strategy, thus we will not enter in details in the description of POD, we only summarize the main aspects: the interested reader can found a complete overview for example in \cite{kunisch,hinze}.

Given a time step $\Delta \tau >0$ (which could be different from $\Delta t$), consider $t_m\in(0,t_f)$ and $\bar{N}$ s.t. $\bar{N} \Delta \tau = t_m$: first the unreduced model (\ref{unreduced}) is solved in $[0,t_m]$, collecting the matrix of snapshots $\mathcal{X}=(\textbf{C}_j)$, where $\textbf{C}_j\in\mathbb{R}^{N_h}$ is the nodal vector of the FE discretization at time $t_j=j\Delta \tau$, $j=0,\ldots,\bar{N}$. After computing the Singular Value Decomposition (SVD) of $\mathcal{X}$, $\mathcal{X} = U S V^t $, a suitable threshold $k$ is chosen. A largely used strategy is to choose $k$ s.t. 
$$\frac{\sum_{i=1}^k S(i,i)^2}{\sum_{i=1}^{min(N_h,\bar{N})} S(i,i)^2}$$ 
is greater than a fixed tollerance. Another possibility is to impose that the first $k$ singular values are greater than a fixed tollerance $\tau_{\sigma}>0$.

It can be proved \cite{hinze} that the \textit{$k$-th POD basis} $\left\{u_i\right\}$, $i=1,\ldots,k$, $u_i:=U(:,i)\in\mathbb{R}^{N_h}$ is the solution of the following minimization problem 

\begin{equation}\min_{\mbox{\boldmath $\xi$}_1\ldots,\mbox{\boldmath $\xi$}_k \in\mathbb{R}^{N_h}}\sum_{j=1}^{\bar{N}} \left|\textbf{C}(t_j)-\sum_{i=1}^k(\textbf{C}(t_j)\cdot \mbox{\boldmath $\xi$}_j)\mbox{\boldmath $\xi$}_i\right|^2,\qquad s.t.\qquad  \mbox{\boldmath $\xi$}_i \cdot \mbox{\boldmath $\xi$}_j=\delta_{ij},\qquad 1\leq i,j\leq k,\label{pr1}\end{equation}
i.e. for every fixed $k$ the mean square error between the elements $\textbf{C}(t_j)$ and the corresponding $k-th$ partial sum of $\sum_{i=1}^k(\textbf{C}(t_j)\cdot \mbox{\boldmath $\xi$}_j)\mbox{\boldmath $\xi$}_i$ is minimized on average.

Finally (\ref{unreduced}) is projected on the space generated by the first $k$ POD basis vectors, i.e. we solve the reduced system
\begin{equation}
\begin{array}{r l l}
U_k^t M U_k \dot{\textbf{a}}(t)+U_k^t A(\textbf{u}) U_k \textbf{a}(t) & = & U_k^t\textbf{F}(C_{in}),\\
\textbf{a}(0)& = & U_k^t\textbf{C}_0.
\label{reduced}
\end{array}
\end{equation}
in $(t_m,t_f)$, where $U_k:= U(:,1:k)$, i.e. the system is projected on the subspace generated by the first $k$ columns of $U$. We denote with $\mbox{\boldmath $\tilde{C}$}(t):=U_k \textbf{a}(t)$ the estimate of $\textbf{C}(t)$, $t\in(t_m,t_f)$ computed using POD.

\section{Continuous inverse problem formulation}
\label{inverso}
We are interested in solving the following inverse problem: given the additional \textit{a priori} information 
\begin{equation}
C = C_s,  \qquad on \quad [0,t_f)\times\Gamma_{down},
\end{equation}
where $C_s\in \mathcal{C}^0((0,t_f),L^2(\Gamma_{down}))$ is a known scalar function, determine $C^*_{in}\in H^{\frac{1}{2}}(\Gamma_{in})$ such that
\begin{equation}
C^*_{in} = \arg \min_{C_{in}\in H^{\frac{1}{2}}(\Gamma_{in})} \mathcal{F}_d(C_{in}),
\label{costo}
\end{equation}
where the \textit{cost function} is
$$\mathcal{F}_d(C_{in}) := \left\|C(C_{in}; t,\textbf{x})-C_s(t,\textbf{x})\right\|^2_{L^2([0,t_f]\times\Gamma_{down})} = \int_0^{t_f} \int_{\Gamma_{down}} (C(C_{in}; t,\textbf{x})-C_s(t,\textbf{x}))^2 d\gamma dt$$
and we have explicited the dependence of $C$, solution of (\ref{direct_problem}), on $C_{in}$: $C(C_{in}; t,\textbf{x}):=C(t,\textbf{x})$ s.t. $C(t,\textbf{x})=C_{in}(\textbf{x})$, when $\textbf{x}\in\Gamma_{in}$ .

As mentioned in \cite{zabaras}, one may consider $C_s$ to be a desired one. In that case, the present inverse problem is a \textit{design problem} where the boundary flux $C_{in}$ is controlled such that a desired concentration is achieved on the boundary $\Gamma_{down}$. $C_s$ can also be considered to represent a continuous approximation of a set of discrete experimental \textit{temperature measurements} obtained at finite number of locations in the boundary $\Gamma_{down}$ and at discrete time instances within the interval $[0,t_f)$. In this paper we mainly refer to the second case. Observe that this class of inverse problems are of significant experimental interest for situations where the direct measurement of the heat flux $C_{in}$ is not possible.

\

As indicated in \cite{zabaras}, the main difficulty with the minimization problem (\ref{costo}) is the calculation of the gradient of $\mathcal{F}$. Mainly two different approaches could be used: the \textit{first discretize than optimize} or vice versa the \textit{first optimize than discretize}. In this paper we focus on the first strategy: in particular we adopt a discrete approximation of $\mathcal{F}^{'}(C_{in})$, combined with a Gauss-Newton approach, as explained starting from section \ref{known}.

\section{Formulation of the discrete inverse problem}
In the \textit{first discretize than optimize} context, $C_s(t,\textbf{x})$ is known only in the $n_y$ nodes of $\Gamma_{down}$, for every discrete time $t_j$, $j=0,\ldots,N-1$. We denote with $\textbf{C}_s(j)\in\mathbb{R}^{n_y}$ the vector of measured concentration at iteration $j$.

As in \cite{girault2008}, we assume that the flow dynamic boundary conditions are steady state and for simplicity we suppose that 
$$\Gamma_{in}=\bigcup_{l=1}^{n_{\theta}}\Gamma^{(l)}_{in},$$ 
being $\Gamma^{(l)}_{in}$ disjoint sets, such that $C_{in}$ is constant on each $\Gamma^{(l)}_{in}$, for all $l=1,\ldots,n_{\theta}$. 

Thus we have to estimate a vector \mbox{\boldmath $\vartheta$} of $n_{\theta}$ non negative parameters: equivalently we assume that the function $C_{in}\in H^{\frac{1}{2}}(\Gamma_{in})$ could be identified by a piecewise constant function $C_{in}(\mbox{\boldmath $\vartheta$})$ such that 
$$C_{in}(\mbox{\boldmath $\vartheta$})(\textbf{x})=\vartheta(l),\qquad \textbf{x}\in\ \Gamma^{(l)}_{in}.$$
Since we are solving an inverse problem we indicate with $\mbox{\boldmath $\hat{\vartheta}$}$ the estimate of the real parameters $\mbox{\boldmath $\vartheta$}$. Let $\Pi:\ \mathbb{R}^{N_h}\rightarrow\mathbb{R}^{n_y}$ be the map which projects the solution of (\ref{unreduced}) on the $n_y$ nodes of $\Gamma_{down}$: thus we denote with $\Pi(\textbf{C}(\mbox{\boldmath $\hat{\vartheta}$};t))$ the estimate of $C(C_{in}(\mbox{\boldmath $\vartheta$}); t,\textbf{x})$ at time $t$ on $\Gamma_{down}$ (\textit{predicted concentration}) obtained solving (\ref{direct_problem}) imposing $C_{in}(\mbox{\boldmath $\hat{\vartheta}$})$ on $\Gamma_{in}$.

In a space-time discrete setting (\ref{costo}) could be restated as
\begin{equation}
\mbox{\boldmath $\vartheta$}^* = \arg \min_{\mbox{\boldmath $\vartheta$}\in \mathbb{R}_+^{n_{\theta}}}
\mathcal{F}_d(\mbox{\boldmath $\vartheta$}),
\label{costo_discreto}
\end{equation}
where the \textit{discrete cost function} is defined as
\begin{equation}\mathcal{F}_d(\mbox{\boldmath $\vartheta$}):= \frac{1}{N}\sum_{j=1}^N \left\|\Pi(\textbf{C}(C_{in}(\mbox{\boldmath $\vartheta$}); j))-\textbf{C}_s(j)\right\|_2^2.\label{costo_d}\end{equation}

\subsection{Proper Orthogonal Decomposition reduction}
Using model order reduction techniques to solve (\ref{costo}), consists in replacing the cost function (\ref{costo_d}) in (\ref{costo_discreto}) with the following one
\begin{equation}
\mathcal{F}_d(\mbox{\boldmath $\vartheta$}):=\frac{1}{N}\sum_{j=1}^N \left\|\Pi(\mbox{\boldmath $\tilde{C}$}(C_{in}(\mbox{\boldmath $\vartheta$}); j))-\textbf{C}_s(j)\right\|_2^2
\label{costo_discreto_ridotto}
\end{equation}
where \mbox{\boldmath $\tilde{C}$} is the solution of (\ref{reduced}). An example of application of POD to solve optimal control problem can be found e.g. in \cite{hinze}.

Since the POD basis depends on the collected snapshots, it is necessary to update the projection space as the estimated control $C_{in}$ varies. Let $\bar{n}$ a small positive integer: at every iteration $i$ in this paper we adopt the following index
$$\mathcal{I}^{(i)}:=\frac{1}{\bar{n}}\left\|\sum_{j=1}^{\bar{n}}\mbox{\boldmath $\tilde{C}$}(C_{in}(\mbox{\boldmath $\vartheta$}^{(i)}); j)-\textbf{C}(C_{in}(\mbox{\boldmath $\vartheta$}^{(i)}); j)\right\|^2_2,$$
i.e. we compare the first iterations of the unreduced system with those obtained projecting on the old POD basis used at iteration $i-1$. Only if $\mathcal{I}^{(i)}$ is greater than a fixed threshold, the $i$-th basis is updated, computing new snapshots, as described in section \ref{sezione_pod}. Two strategies can be used \cite{hinze}: old snapshots can be discarded or not. In practice this consists in adding POD modes computed in the $i-1$-th iteration to the new snapshots ensamble: in this case the projection space is more robust to control variations but usually is slightly bigger. For our experimental tests we prefer to discard old snapshots. We obserse that in \cite{hinze} a new basis is computed at every iteration, without considering an index $\mathcal{I}$.

\section{Known source location $\Gamma_{in}$}
\label{known}
As a first step toward the solution strategy, we consider a simpler problem, assuming that the source location $\Gamma_{in}$ is known.  

\subsection{Solution uniqueness}
In this section we demonstrate that if $\Gamma_{in}$ is known, then the discrete inverse problem admits a unique solution, since there are no local minima. Moreover changes in $C_{in}$ corresponds to changes in the registered concentration.

First of all we prove the following Lemma, which justifies mathematically the physical principle that, as $C_{in}$ increases on $\Gamma_{in}$, the concentration on $\Gamma_{down}$ increases too.

\begin{figure}[h]
\begin{center}
\includegraphics[width=10cm]{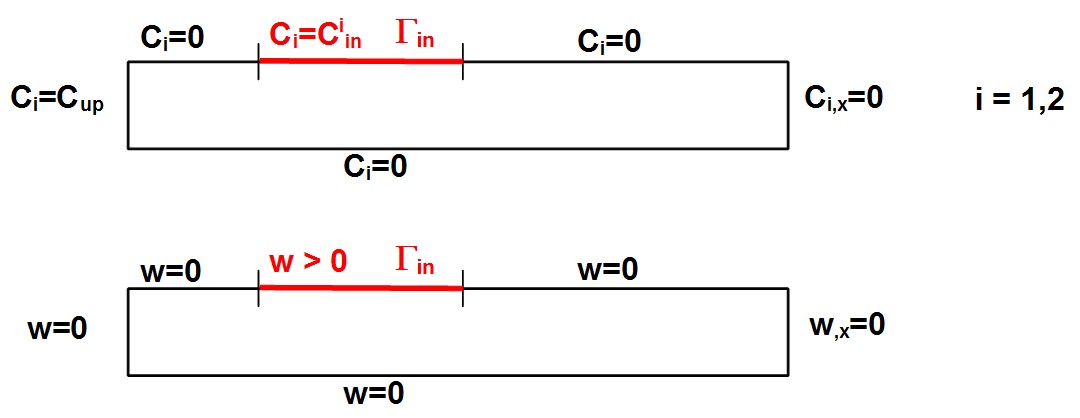}
\end{center}
\caption{\small \textit{}}\normalsize
\label{studio_domini}\end{figure}

\begin{LEM}
\label{monotonia}
Consider the two problems
\begin{equation}
\left\{\begin{array}{r l l l}
\frac{\partial C_i}{\partial t}-\mu \Delta C_i + \nabla \cdot (\textbf{u} C_i)+\sigma C_i & = & 0, & \qquad in \quad (0,t_f)\times\Omega \\ [5 pt]
C_i & = & C_0, & \qquad on \quad \left\{0\right\}\times \Omega \\ 
C_i & = & C^{(i)}_{in}, & \qquad on \quad (0,t_f)\times\Gamma_{in} \\ 
C_i & = & C_{up}, & \qquad on \quad (0,t_f)\times\Gamma_{up} \\ 
\mu\frac{\partial C_i}{\partial n} & = & 0, & \qquad on \quad (0,t_f)\times\Gamma_{down}\\ 
C_i & = & 0, & \qquad on \quad (0,t_f)\times\Gamma_{r}
\end{array}\right.
\end{equation}
represented in Figure \ref{studio_domini} (up), where $i=1,2$. Suppose that $C^{(2)}_{in}(\textbf{x})>C^{(1)}_{in}(\textbf{x})$, for every $\textbf{x}\in\Gamma_{in}$. 

Then $C_2(t,\textbf{x})>C_1(t,\textbf{x})$ for every $t\in(0,t_f)$ and $\textbf{x}\in \Gamma_{down}$.
\end{LEM}

\textbf{Proof.}

Define $w:=C_2-C_1$, which solves 
\begin{equation}
\left\{\begin{array}{r l l l}
\frac{\partial w}{\partial t}-\mu \Delta w + \nabla \cdot (\textbf{u} w)+\sigma w & = & 0, & \qquad in \quad (0,t_f)\times\Omega \\ [5 pt]
w & = & 0, & \qquad on \quad \left\{0\right\}\times\Omega \\ 
w & = & C^{(2)}_{in}-C^{(1)}_{in}, & \qquad on \quad (0,t_f)\times\Gamma_{in} \\ 
w & = & 0, & \qquad on \quad (0,t_f)\times\Gamma_{up} \\ 
\mu\frac{\partial w}{\partial n} & = & 0, & \qquad on \quad (0,t_f)\times\Gamma_{down}\\ 
w & = & 0, & \qquad on \quad (0,t_f)\times\Gamma_{r}
\end{array}\right.
\label{direct_problem_2}
\end{equation}
as illustrated in figure \ref{studio_domini} (down). 

Observe that $w$ is smooth only inside the domain, but it is not continuous near the boundary, where it admits discontinuities of the first kind: thus generalized solutions must be considered. The strong minimum principle for parabolic operators can be extended for generalized solutions \cite{friedman,ladyzenskaja} (cfr. appendix \ref{max_op}): thus we know that the minimum is assumed at the boundary.

Moreover, for every open neighbourhood $U$ of $\Gamma_{down}$, such that $w$ is regular inside $U \cap \Omega$, $\frac{\partial w}{\partial n}= 0$ on $\Gamma_{down}$ implies that the maximum and the minimum of $w$ over $U \cap \Omega$ cannot belong to $\Gamma_{down}$ (cfr. \cite{friedman}). 

As a consequence $w\geq 0$ in $(0,t_f)\times \Omega$ and, since the minimum is not attained on $\Gamma_{down}$, $w=C_2-C_1>0$ on $\Gamma_{down}$, for all $t\in(0,t_f)$ i.e. the thesis holds true. 

\begin{flushright}
$\square$
\end{flushright}

The following Proposition is equivalent to prove that there are no local minima.
\begin{PROP}\label{convessita_cd}
For every $\mbox{\boldmath $\bar{\vartheta}$}\in \mathbb{R}_+^{n_{\theta}}$, $\mbox{\boldmath $\bar{\vartheta}$}\neq \mbox{\boldmath $\vartheta^{*}$}$, there exists at least a sequence of profiles $\left\{\mbox{\boldmath $\vartheta$}\right\}_n$, $\mbox{\boldmath $\vartheta$}_0=\mbox{\boldmath $\bar{\vartheta}$}$, converging in $\mathcal{L}^2(\mathbb{R}^{n_{\theta}})$ to the real profile $\mbox{\boldmath $\vartheta$}^{*}$, such that $\mathcal{F}_d(\mbox{\boldmath $\vartheta$}_n)\downarrow \mathcal{F}_d(\mbox{\boldmath $\vartheta$}^*)$.
\end{PROP}

\textbf{Proof.} We can construct the sequence $\left\{\mbox{\boldmath $\vartheta$}_n\right\}_n$ in the following way. For every $k= 1,\ldots,n_{\theta}$ 
\begin{equation}\vartheta_k(j):=\left\{\begin{array}{l l} 
\vartheta_{k-1}(j),& j\neq k\\
 \bar{\vartheta}(j) - (\bar{\vartheta}(j)- \vartheta^*(j)),& j=k \end{array}\right..\label{successione}\end{equation}
Thus $\mbox{\boldmath $\vartheta$}_{n_{\theta}}=\mbox{\boldmath $\vartheta$}^*$ by construction. Moreover the corresponding sequence of cost functions is \textit{decreasing}: $\mathcal{F}_d(\mbox{\boldmath $\vartheta$}_1)>\mathcal{F}_d(\mbox{\boldmath $\vartheta$}_2)>\ldots>\mathcal{F}_d(\mbox{\boldmath $\vartheta$}^*)$. This fact is a direct consequence of the application of Lemma \ref{monotonia}: suppose that $\vartheta_{k-1}(k)<\vartheta^*(k)$. Then $\vartheta_{k}(k)>\vartheta_{k-1}(k)$ by construction and thus $\Pi(\textbf{C}(\mbox{\boldmath $\vartheta$}_k;t))$ will be higher than $\Pi(\textbf{C}(\mbox{\boldmath $\vartheta$}_{k-1};t))$ for every $t\in(0,t_f)$ (Lemma \ref{monotonia}) and thus closer to $\Pi(\textbf{C}(\mbox{\boldmath $\vartheta$}^*;t))$. Analogously if $\vartheta_{k-1}(k)>\vartheta^*(k)$, applying Lemma \ref{monotonia}, $\Pi(\textbf{C}(\mbox{\boldmath $\vartheta$}_k;t))$ will be lower than $\Pi(\textbf{C}(\mbox{\boldmath $\vartheta$}_{k-1};t))$ for all $t$ and thus closer to $\Pi(\textbf{C}(\mbox{\boldmath $\vartheta$}^*;t))$.

\begin{flushright}
$\square$
\end{flushright}

\subsection{Numerical solution strategy}

Starting from an initial guess $\mbox{\boldmath $\hat{\vartheta}$}^{(0)}$, line search algorithms find the $k+1$-iteration starting from the $k$-th one in the following way:
$$\mbox{\boldmath $\hat{\vartheta}$}^{(k+1)} = \mbox{\boldmath $\hat{\vartheta}$}^{(k)}+\alpha^{(k)}\textbf{s}^{(k)},$$
where the \textit{damping parameter} $\alpha^{(k)}$ is obtained using a bisection procedure.

The standard Newton step consists in solving at each iteration the system
$$\mathcal{F}_d^{''}(\mbox{\boldmath $\hat{\vartheta}$}^{(k)}) \textbf{s}^{(k)} = - \mathcal{F}_d^{'}(\mbox{\boldmath $\hat{\vartheta}$}^{(k)}).$$
Let $\mathcal{R}:\ \mathbb{R}^{n_y \times N} \rightarrow\mathbb{R}^{n_y N}$ be the map s.t. starting from an $n_y \times N$ matrix $B=[b_1,\ldots,b_N]$, it gives $\mathcal{R}(B):= \left(\begin{array}{c}b_1\\\vdots\\b_N\end{array}\right)$.
  
The computation of $\mathcal{F}_d^{''}$, Hessian of the cost function, usually is expensive. Moreover, since we are dealing with a \textit{least squares} problem (\ref{costo}), we adopt the Gauss-Newton approximation (cfr. \cite{nocedal}), i.e. we solve
\begin{equation}\Psi_{\mbox{\boldmath $\hat{\vartheta}$}^{(k)}} \textbf{s}^{(k)} =  \textbf{e}_{\mbox{\boldmath $\hat{\vartheta}$}^{(k)}},\label{system_eq}\end{equation}
where the \textit{sensitivity matrix} $\Psi_{\mbox{\boldmath $\hat{\vartheta}$}^{(k)}}\in\mathbb{R}^{n_y N\times n_{\theta}}$ is such that 
$$\Psi_{\mbox{\boldmath $\hat{\vartheta}$}^{(k)}}(:,i):=\frac{\partial}{\partial \mbox{\boldmath $\hat{\vartheta}$}^{(k)}(i)} \mathcal{R}(\Pi(\textbf{C}(\mbox{\boldmath $\hat{\vartheta}$}^{(k)};\cdot))),$$
for all $i=1,\ldots,n_{\theta}$ and the \textit{prediction error} 
$$\textbf{e}_{\mbox{\boldmath $\hat{\vartheta}$}^{(k)}} := \mathcal{R}(\textbf{C}_s(\cdot)) - \mathcal{R}(\Pi(\textbf{C}(\mbox{\boldmath $\hat{\vartheta}$}^{(k)};\cdot))).$$

To compute numerically the sensitivity matrix a finite difference scheme is needed: 
$$\Psi_{\mbox{\boldmath $\hat{\vartheta}$}^{(k)}}(:,j) \approx \frac{1}{\delta}\left[\mathcal{R}(\Pi(\textbf{C}(\mbox{\boldmath $\hat{\vartheta}$}^{(k)}(1),\ldots,\mbox{\boldmath $\hat{\vartheta}$}^{(k)}(j)+ \delta,\ldots,\mbox{\boldmath $\hat{\vartheta}$}^{(k)}(n_{\theta});\cdot)))- \mathcal{R}(\Pi(\textbf{C}(\mbox{\boldmath $\hat{\vartheta}$}^{(k)};\cdot)))\right],$$
where $\delta>0$ is a small perturbation parameter.

Observe that in general this approximation is computationally expensive, since, it requires the computation of the concentration also for the perturbed input. When $\Gamma_{in}$ is known, only very few parameters are considered, thus this approximation is effective. The problem becomes more involving when $\Gamma_{in}$ is unknown, since the number of parameters is higher: in section \ref{pa} we will explain how the adaptive parametrization and time localization can reduce the computational cost. 

If $\delta>0$ is too small, the finite difference estimate could be inaccurate, since at the numerator we are considering the difference between two quantities which has approximately the same absolute value, and this is divided by a very small denominator, which amplifies the error. A possible solution e.g. is to adopt the Complex-Step Derivative Approximation \cite{martins}, in which an immaginary increment $i\delta$ is used, approximating $$\Psi_{\mbox{\boldmath $\hat{\vartheta}$}^{(k)}}(:,j)\approx \frac{1}{\delta} Im \left(\mathcal{R}(\Pi(\textbf{C}(\mbox{\boldmath $\hat{\vartheta}$}^{(k)}(1),\ldots,\mbox{\boldmath $\hat{\vartheta}$}^{(k)}(j)+ i\delta,\ldots,\mbox{\boldmath $\hat{\vartheta}$}^{(k)}(n_{\theta});\cdot)))\right).$$ 

Finally observe that we are assuming that the pollutant is put into the domain, thus 
$$C_{in}\geq 0:$$
as a consequence we need also a \textit{projection} step onto $[0,+\infty)$ of each component of $\mbox{\boldmath $\hat{\vartheta}$}^{(k)}$, after its computation. 

\subsection{Numerical results}
\label{numerical_results_known}
In this section the Projected damped Gauss Newton is compared to other classical solution strategies. Experimental data are simulated numerically, on $\Omega = [0,8]\times [0,1]$,  $\Gamma_h=[0,8]\times\left\{1\right\}\cup [0,8]\times\left\{0\right\}$. Moreover the velocity field $\textbf{u}$ is modelled as a Poiseuille flow i.e.
$$\textbf{u}(x_1,x_2)=\left(\begin{array}{c c} -4\nu x^2_2 + 4 \nu x_2\\ 0 \end{array}\right).$$
We assume that $\nu=50$, $\mu=0.1$, $\sigma=0.1$ and $C_{up}=0.1$. Moreover in this section a Gaussian error of variance $0.05$ and mean zero is added. 

Classical solution strategies cited in this section are well described e.g. in \cite{kaltenbacher}. As a regularization parameter, when needed, we use $\alpha=0.01$, moreover we choose a maximum number of iterations $max_{it}=20$. Consider the following two \textit{sparse} examples:
\begin{enumerate}
	\item $\Gamma_{in}=[4,4.5]\times\left\{1\right\}$, $\vartheta = 100$; \\
	\item $\Gamma_{in}=[4.5,5]\times\left\{1\right\}\cup[1.5,2]\times\left\{0\right\},\ \mbox{\boldmath $\vartheta$} = (100,80)$;\\
\end{enumerate}
and see how different techniques approximate them.

\begin{figure}[t]
\begin{center}
\includegraphics*[width=3.5cm]{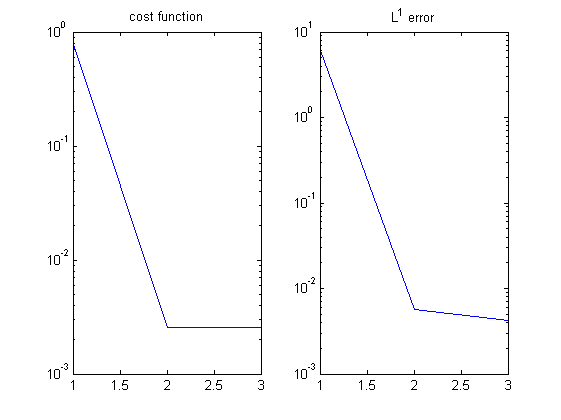}
\includegraphics*[width=3.5cm]{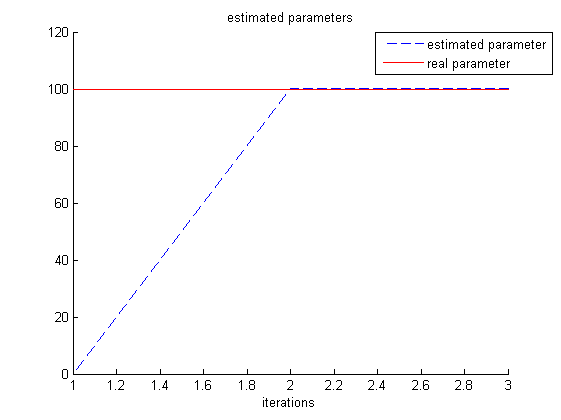}

\includegraphics*[width=3.5cm]{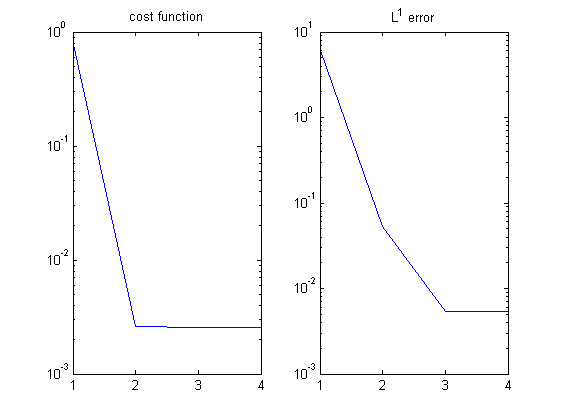}
\includegraphics*[width=3.5cm]{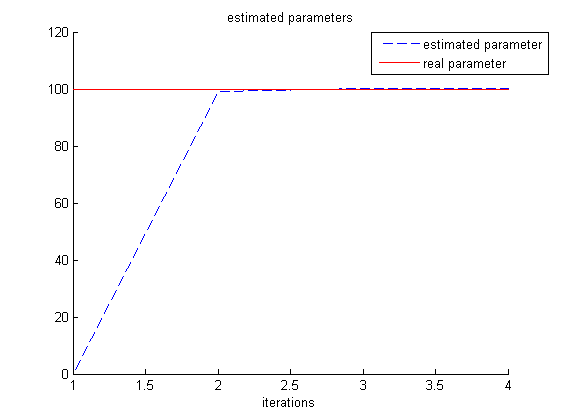}

\includegraphics*[width=3.5cm]{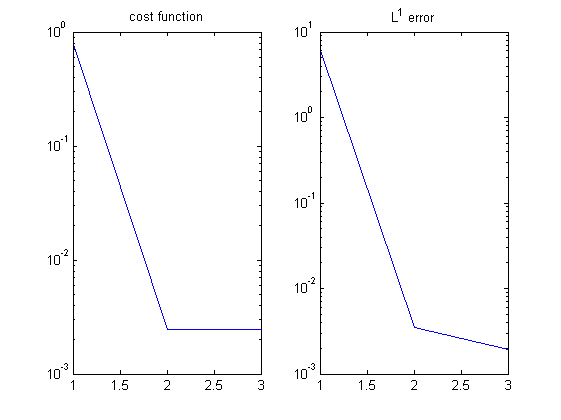}
\includegraphics*[width=3.5cm]{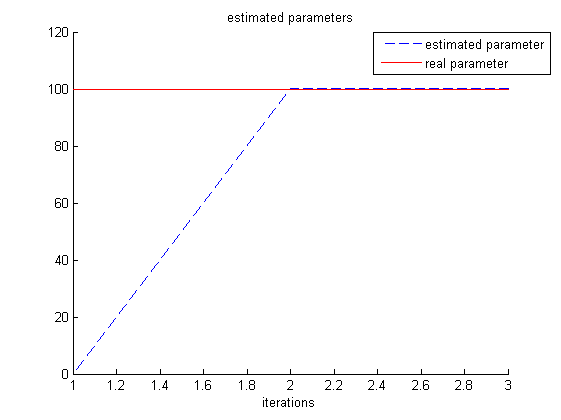}

\includegraphics*[width=3.5cm]{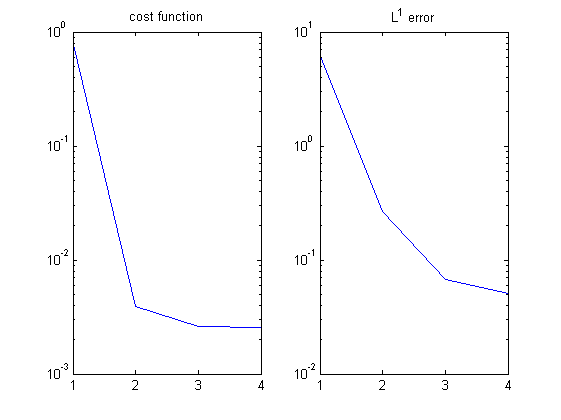}
\includegraphics*[width=3.5cm]{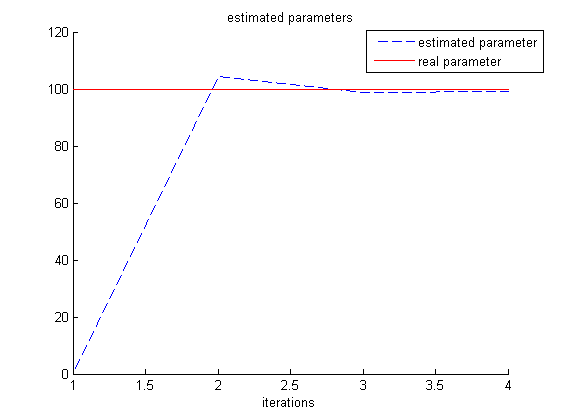}
\end{center}
\caption{\small \textit{First example: different strategies. Left: cost function and error, right: convergence. First row: projected damped Gauss Newton, second row: Levenberg Marquardt, third row: steepest descent, fourth row: Tikhonov method.}}\normalsize
\label{1}
\end{figure}

\begin{figure}[t]
\begin{center}
\includegraphics*[width=3.5cm]{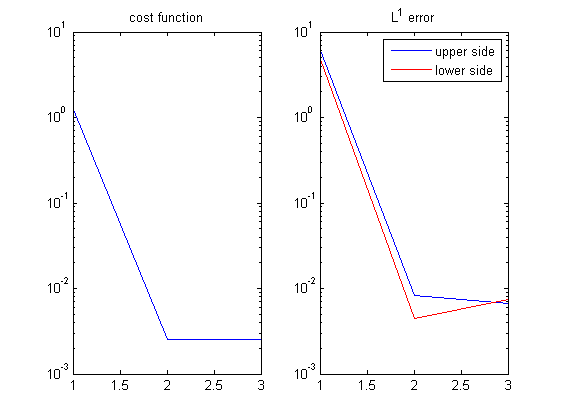}
\includegraphics*[width=3.5cm]{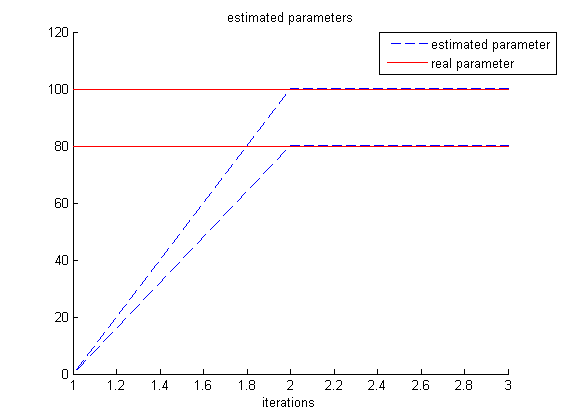}

\includegraphics*[width=3.5cm]{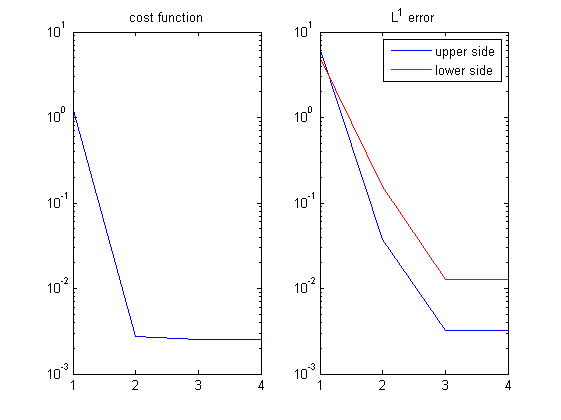}
\includegraphics*[width=3.5cm]{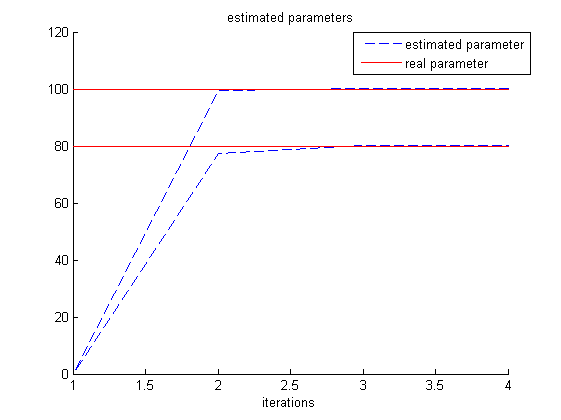}

\includegraphics*[width=3.5cm]{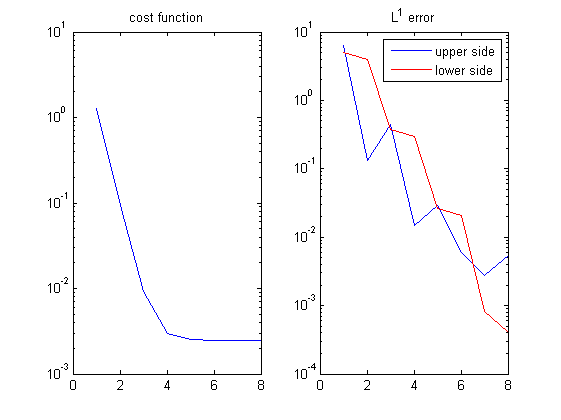}
\includegraphics*[width=3.5cm]{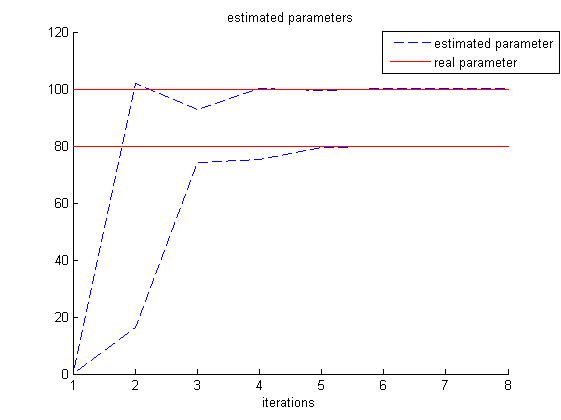}

\includegraphics*[width=3.5cm]{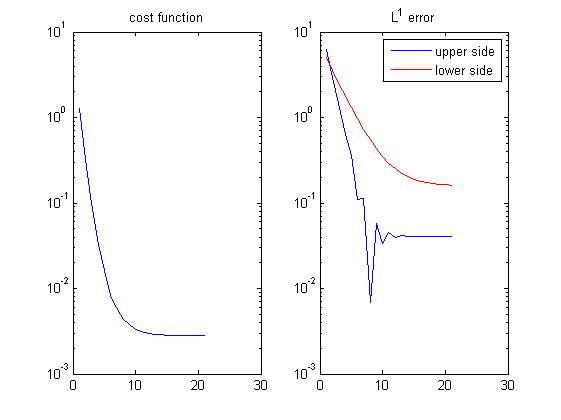}
\includegraphics*[width=3.5cm]{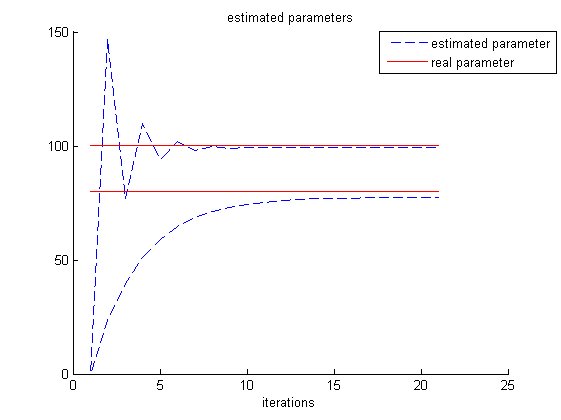}
\end{center}
\caption{\small \textit{Second example: different strategies. Left: cost function and error, right: convergence. First row: projected damped Gauss Newton, second row: Levenberg Marquardt, third row: steepest descent, fourth row: Tikhonov method.}}\normalsize
\label{2}
\end{figure}

First of all we consider the example 1. Performances of different methods are depicted in figure \ref{1}. 
In the second example, two parameters have to be estimated: results are plotted in figure \ref{2}.

Observe that in both cases the Projected damped Gauss Newton algorithm performs well, converging faster to the optimal solution. It should be noted that, in constrast to Tikhonov and Levenberg-Marquardt it does not need a regularization parameter: this is important because it tells us that, knowing the source location, the inverse problem is not ill-posed, as stated in Proposition \ref{convessita_cd}. 

\subsection{Reduce the order of the system using POD}
In this section we analyze the POD reduction introduced in section \ref{sezione_pod} on a test case. Consider example 2 introduced in the previous section; in POD reduction two parameters plays a central role: $t_m$, which characterizes the interval $[t_0,t_m]$ when snapshots are collected, and the threshold $\tau_{\sigma}$ on the singular values of the snapshots matrix. As can be seen in table \ref{pod_ese2}, increasing $t_m$ corresponds to a better approximation, since more snapshots are collected. To obtain higher accuracy decreasing $t_m$, it is necessary to increase $\tau_{\sigma}$, i.e. to consider a bigger reduced model. 
\begin{center}
\begin{table}[h]
\scriptsize{\begin{tabular}{|c|c||c|c|c|c|c|} \hline
\multicolumn{2}{|c||}{}&\multicolumn{2}{|c|}{$L^1$ error:} & \multicolumn{1}{|c|}{$\mathcal{F}_d(\mbox{\boldmath $\vartheta$})$} & \multicolumn{1}{|c|}{Dim. model} &  \multicolumn{1}{|c|}{num. it.}\\ 
\multicolumn{2}{|c||}{} & \multicolumn{1}{|c|}{up} & \multicolumn{1}{|c|}{down} & \multicolumn{1}{|c|}{} & \multicolumn{1}{|c|}{} &  \multicolumn{1}{|c|}{} \\ \hline
\multicolumn{2}{|c||}{Unreduce model} & 0 & 0 & $10^{-20}$ & 1071 & 2 \\ \hline\hline
\multicolumn{2}{|c||}{Reduced models:} & \multicolumn{5}{|c|}{}\\
\multicolumn{1}{|c|}{$t_m$} & \multicolumn{1}{|c||}{$\tau_{\sigma}$} & \multicolumn{5}{|c|}{}\\
\hline
2.5 & 0.01 & 0.117 & 4.8 & 0.113 & 23 & 5 \\ \hline
2.5 & $10^{-4}$ & 0.083 & 1.24 & 0.089 & 32 & 4 \\ \hline
3.75 & 0.01 & 0.08 & 0.08 & $6\cdot 10^{-4}$ & 25 & 3 \\ \hline
3.75 & $10^{-4}$ & 0.02 & 0.02 & $3\cdot 10^{-5}$ & 39 & 4 \\ \hline
5 & 0.01 & 0.0015 & 0.0015 & $10^{-6}$ & 29 & 3\\ \hline
\end{tabular}}
\caption{\small \textit{Example 2 of section \ref{numerical_results_known}, choosing different intervals $[t_0,t_m]$ to collect snapshots and different thresholds $\tau_{\sigma}$ on singular values of the snapshots matrix.}}  
\label{pod_ese2}
\end{table}
\end{center}
It is important to note that the reduction is significative with respect to the unreduced model, which has dimension $1071$. However, as described in section \ref{sezione_pod}, it should be noted that it is necessary to update the POD basis: in all these examples the basis is updated at every  new iteration, imposing $0.1$ as a threshold on $\mathcal{I}^{(i)}$.

\

A more involving problem is considered in section \ref{pa}, where it is assumed that also the source location $\Gamma_{in}$ is unknown. In general in that case projected damped Gauss Newton could not be sufficient and it is too costly, thus it is necessary to adopt a suitable solution strategy based upon an adaptive parametrization and time localization.

\section{Unknown source location $\Gamma_{in}$}
\label{pa}
Suppose now that the \textit{location} of $\Gamma_{in}$ is unknown.

\subsection{Ill-posedness of the problem}
\label{fiume_1d}
To study analytically what happens when $\Gamma_{in}$ is unknown, we consider a simplified model problem: let $C=C(x)$, $x\in [x_1,x_2]\subset\mathbb{R}$, $x_2>x_1$ be the solution of the following one dimensional ODE:
\begin{equation}
\left\{\begin{array} {r l l}
-\mu C^{''}(x) + u C^{'}(x) &=& f(x),\qquad in\ (x_1,x_2), \\
C(x_1) &=& C_{up},\\
C^{'}(x_2)&=&0,
\end{array}\right.
\label{ode1d}
\end{equation}
where $ f(x) = \left\{\begin{array}{r l} M, & \left|x-x_m\right|\leq h\\
0,& elsewhere\ in\ (x_1,x_2) \end{array}\right.$, $M>0$, $x_m\in (x_1,x_2)$, $h \in (0,1)$ s.t. $x_m \pm h \in (x_1,x_2)$.
Observe that (\ref{ode1d}) can be viewed as the one dimensional stationary counterpart of (\ref{direct_problem}) when $\sigma=0$ and considering only the $x$-axis in figure \ref{dominio_fiume}: the unknown immision boundary $\Gamma_{in}$ can be represented by an unknown forcing term $f$, applied in $[x_m-h, x_m+h]$, of intensity $M$. In this context the inverse problem (\ref{costo}) is equivalent to determine the source position ($h$ and $x_m$) and intensity ($M$) given the measured concentration $C_s\in\mathbb{R}$ in $x_2$.

Problem (\ref{ode1d}) can be solved analytically, obtaining 
$$ C(x) = \left\{\begin{array}{r l l} 
c_1 + c_2 e^{\frac{u}{\mu}x}, & x < x_m-h,\\
c_3 + \frac{M}{u} x + c_4 e^{\frac{u}{\mu}x}, & \left|x-x_m\right|\leq h,\\
c_5 + c_6 e^{\frac{u}{\mu}x}, & x > x_m+h,
\end{array}\right.$$
where $c_1,\ldots,c_6$ are suitable real coefficients obtained imposing boundary conditions and continuity of $u$ and $u^{'}$ in $x_m\pm h$. In particular we are interested in estimating the concentration at the measurement point $x=x_2$. For simplicity we assume that $x_1=0$ and $x_2=1$. Thus it can be derived that
$$c_6 = 0,\qquad c_5 = \frac{1}{u^2}\exp{\frac{-u(x_m+h)}{\mu}}\left(2 u h M \exp{\frac{u(x_m+h)}{\mu}}+\mu M \left(1-\exp{\frac{2 u h}{\mu}}\right)\right).$$ 
\begin{figure}[h]
\begin{center}
\includegraphics[width=3.5cm]{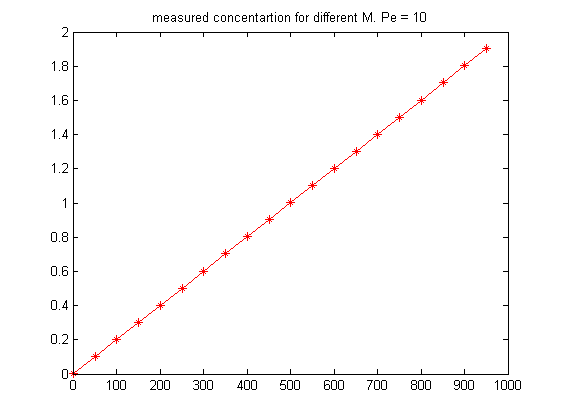}
\includegraphics[width=3.5cm]{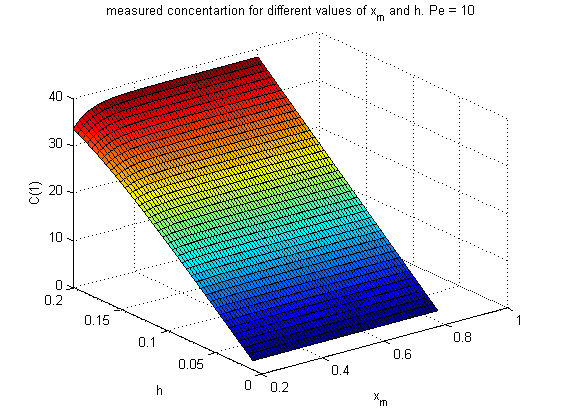}
\includegraphics[width=3.5cm]{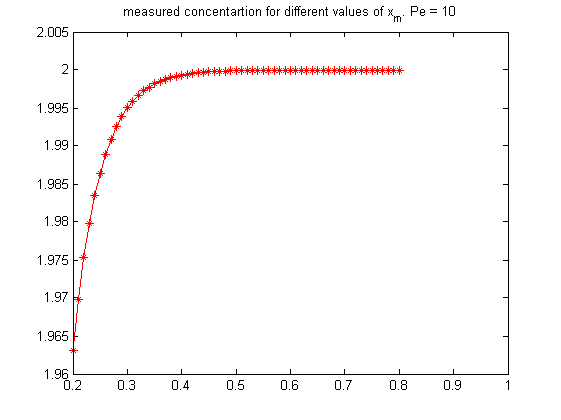}
\end{center}
\caption{\small \textit{Solution of (\ref{ode1d}) at the measurement point $x_2=1$ for different values of $M$ (left), $h$ and $x_m$ (center), $x_m$ (right).}}\normalsize
\label{figure_ode}\end{figure}
Thus $C(x)$ is constantly equal to $c_5$ in $[x_m+h,1]$. We can now study how $C(1)$ depends on $M$, $h$, $L$. We consider $\mu=0.5$ and $u=10$ (Peclet number $Pe=\frac{u}{2\mu}=10$, quantity that characterize convection diffusion problems). As can be seen in figure \ref{figure_ode}, varying only $M$, fixing $h$ and $x_m$ (i.e. knowing the source location), corresponds to a linear striclty increasing $C(1)$ (cfr. figure \ref{figure_ode} (left)). On the contrary fixing $M$ but varying $h$ and $x_m$ corresponds to the surface plotted in figure \ref{figure_ode} (center): fixing $h$ for different values of $x_m$ we obtain almost the same $C(1)$ (cfr. figure \ref{figure_ode} (right)). Thus measuring $C(1)$, \textit{the problem of determining the source is ill-posed in the stationary regime}. Increasing the Peclet number this phenomenum is stressed.

Even for this simplified 1D stationary problem, in general unknown source position gives rise to an ill-conditioned problem. 

\subsection{Numerical solution of the discrete inverse problem}

The problem consists in estimating both \textit{the position} of the sources $C_{in}$ in the horizontal segments $\Gamma_h:=\Gamma_r \cup \Gamma_{in}$ and their \textit{intensity}. 

\subsubsection{Algorithm 1: working on the finest subdivision}
First of all we consider $\left\{x_1,\ldots,x_{\frac{n_{\theta}^{(f)}}{2}+1}\right\}$ a \textit{reference uniform finest subdivision} of $\Gamma_{h}$ of step length $\Delta x$, which represents the minimum amplitude of estimated source emissions. 

The simplest strategy consists in applying the Gauss Newton method directly on the finest subdivision (cfr. algorithm \ref{finest}), i.e. in estimating $n_{\theta}^{(f)}$ parameters. This problem is particularly demanding for its high computational cost, due to the large number of parameters to be estimated at each Newton's iteration. Moreover it should be noted that if we want to estimate a \textit{sparse} vector of parameters, working only on the finest subdivision is not efficient, as we will see in the following sections.
\begin{algorithm}
	\footnotesize
	\caption{\small Sketch of the algorithm working on the finest subdivision:}
\label{finest}
  \begin{algorithmic}[1]
\STATE{Given the finest subdivision of $\Gamma_h$, $\mbox{\boldmath $\hat{\theta}$}^0=\textbf{0}$, $\mu^0=1$;}  
\WHILE{$\mathcal{F}_d(\mbox{\boldmath $\hat{\vartheta}$}^{l}) < tol$}  
\STATE {solve $\psi_{\hat{\theta}^k} \textbf{s}^k =  \textbf{e}_{\hat{\theta}^k}$ ;}
\STATE{$\mbox{\boldmath $\hat{\theta}$}^{k+1}=\mbox{\boldmath $\hat{\theta}$}^k + \mu^{k} \textbf{s}^k$}
\STATE{\textit{projection}: for every $j\in[0,n_{\theta}-1]$ s.t. $\hat{\theta}^{k+1}(j)<0$, impose $\hat{\theta}^{k+1}(j)=0$}
\STATE{compute $\mathcal{F}_n(\mbox{\boldmath $\hat{\theta}$}^{k+1})$}
\IF{$\mathcal{F}_n(\mbox{\boldmath $\hat{\theta}$}^{k+1})>\mathcal{F}_n(\mbox{\boldmath $\hat{\theta}$}^{k})$}
\STATE{$l=0$;}
\STATE{$\mu^{k,l}=\frac{\mu^k}{2}$}
\WHILE{$\mathcal{F}_n(\mbox{\boldmath $\hat{\theta}$}^{k+1})<\mathcal{F}_n(\mbox{\boldmath $\hat{\theta}$}^{k})$}
\STATE{$\mbox{\boldmath $\hat{\theta}$}^{k+1}=\mbox{\boldmath $\hat{\theta}$}^k + \mu^{k,l} \textbf{s}^k$}
\STATE{$l=l+1$;}
\STATE{$\mu^{k,l}=\frac{\mu^{k,l}}{2}$}
\ENDWHILE
\ENDIF
\ENDWHILE
\end{algorithmic}		
\end{algorithm}
To solve the system (\ref{system_eq}) both \textit{TSVD} and \textit{diagonal scaling} \cite{nocedal} are used. The last one, presented in \cite{marcuzzi} to solve a conduction inverse problem, works as follows: at iteration $k$, given the subdivision $\mathcal{S}^{(k)} = \left\{x_1,\ldots,x_{\frac{n_{\theta}^{(f)}}{2}+1}\right\}$, for every $i=1,\ldots,n^{(k)}_c$, where $n^{(k)}_c=n_{\theta}^{(f)}$ denotes the number of columns of $\Psi_{\mbox{\boldmath $\hat{\vartheta}$}^{(k)}}$ at iteration $k$, $\Psi_{\mbox{\boldmath $\hat{\vartheta}$}^{(k)}}(:,i)$ is multiplied by a weight $d_i$, equal to the length of the maximal segment of the current subdivision, divided by the length of the segment corresponding to the $i$-th column. Thus diagonal scaling corresponds to solve
\begin{equation}\begin{array}{r l l}
\Psi_{\mbox{\boldmath $\hat{\vartheta}$}^{(k)}} D^{(k)} \mbox{\boldmath $\tilde{s}$}^{(k)} & = &  \textbf{e}_{\mbox{\boldmath $\hat{\vartheta}$}^{(k)}},\qquad D^{(k)}=diag(d_i^{(k)}),\qquad d_i^{(k)}=\frac{\max_{x^{k}_{j+1},x^{k}_{j}\in\mathcal{S}^{(k)}} x^{k}_{j+1}-x^{k}_{j}}{x^{k}_{i+1}-x^{k}_{i}},\\
\textbf{s}^{(k)} &=& D^{(k)} \mbox{\boldmath $\tilde{s}$}^{(k)},
\end{array}\end{equation}
instead of (\ref{system_eq}).

\subsubsection{Algorithm 2: working on the finest subdivision with time localization}
As explained in section \ref{fiume_1d}, in the stationary regime the problem is illposed: time localization corresponds to a better conditioned problem, since it consists in selecting only those rows of the sensitivity matrix which are significative for the dynamic, i.e. corresponding to the transitional dynamics.

\begin{figure}[h]
\begin{center}
\includegraphics*[width=10cm]{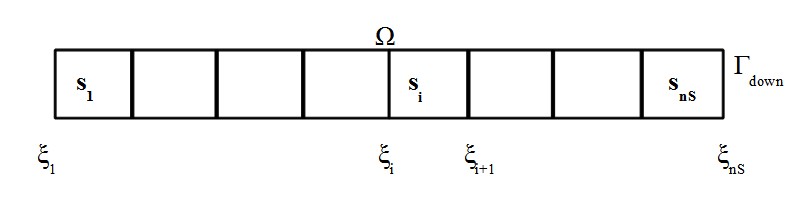}
\end{center}
\caption{\small \textit{Example of partition of $\Omega$ in sections.}}\normalsize
\label{esempio_sezioni_cd}
\end{figure}
More precisely, the idea is to partition the domain $\Omega$ in a suitable number $n_s>1$ of \textit{sections} $\mathcal{U}=\left\{s_j\right\}$, $j=1,\ldots,n_s$ (cfr. e.g. figure \ref{esempio_sezioni_cd}). Referring to figure \ref{dominio_fiume}, we suppose that $s_j:=[\xi_j,\xi_{j+1}]\times[y_1,y_2]$, $\xi_1=x_1$, $\xi_{n_s+1}=x_2$. In particular in algorithm 2 we assume that $\left\{\xi_1,\ldots,\xi_{n_s+1}\right\}=\left\{x_1,\ldots,x_{\frac{n_{\theta}^{(f)}}{2}+1}\right\}$, i.e. it coincides with the finest subdivision. Denote with $I^{(j)}$, the set of parameters belonging to $s_{j}$.

The algorithm works as follows: starting from $s_{n_s}$, it computes the sensitivity matrix only of those parameters belonging to $I^{(n_s)}$, only in the time interval $[t_{0}^{(n_s)},t_{f}^{(n_s)}]$, $t_{0}^{(n_s)}\geq t_0$, $t_{f}^{(n_s)}\leq t_f$: below it is explained how to choose the interval. The estimate of the parameters of section $n_s$ is done as explained before, using a projected damped Gauss-Newton method. To regularize the problem both TSVD and diagonal scaling are used. 

Define $\mathcal{O}^{(j)}$, $j=1,\ldots,n_s-1$, the set of parameters estimated in section $s_{j+1}$ greater than a threshold $\epsilon_3>0$: then in section $s_j$ all parameters belonging to $\mathcal{O}^{(j)}\cup I^{(j)}$ will be estimated, only in the time interval $[t_{0}^{(j)},t_{f}^{(j)}]$. In algorithm \ref{time_loc} previous ideas are summarized.

\begin{algorithm}
	\footnotesize
	\caption{\small Sketch of the algorithm working on the finest subdivision with time localization:}
\label{time_loc}
  \begin{algorithmic}[1]
\STATE{Given $\left\{\xi_1,\ldots,\xi_{n_s+1}\right\}$ coincident with the finest subdivision of $\Gamma_h$ and the threshold $\epsilon_3>0$;}  
\WHILE{$\mathcal{F}_d(\mbox{\boldmath $\hat{\vartheta}$}^{l}) < tol$}  
\STATE{$i=n_s$; $\mathcal{O}^{(n_s)}=\emptyset$}
\WHILE{$i>0$}
\STATE {Let $I^{(i)}$ be the set of parameters of $\mbox{\boldmath $\hat{\vartheta}$}^{l}$ that belongs to section $i$;}
\STATE{in $[t_{0}^{(i)},t_{f}^{(i)}]$ apply the regularized projected damped Gauss Newton method to optimize parameters whose indices belong to $I^{(i)} \cup \mathcal{O}^{(i)}$;}
\STATE{update the positions $I^{(i)} \cup \mathcal{O}^{(i)}$ of $\mbox{\boldmath $\hat{\vartheta}$}^{l}$;}
\STATE{define $\mathcal{O}^{(i-1)}$ as the set of indices of parameters greater than $\epsilon_3$;}
\STATE{$i=i-1$;}
\ENDWHILE
\ENDWHILE
\end{algorithmic}		
\end{algorithm}

\subsubsection{Algorithm 3: using an adaptive parametrization}
The alternative is to use an \textit{adaptive parametrization}, i.e. to adaptively update the subdivision of $\Gamma_h$ used in the current iteration of the Newton method. This strategy is particularly indicated when dealing with a \textit{sparse} vector of parameters: in this situation it is important to localize $\Gamma_{in}$ in $\Gamma_h$ and to refine the parametrization possibly only around that point. It reduces the computational cost reducing the number of columns of the sensitivity matrix. A similar strategy has been presented in \cite{deolmi}, to solve an inverse conduction problem of corrosion estimation.  

The algorithm works as follows: starting from an initial \textit{coarse} subdivision of $\Gamma_h$, $\mathcal{S}^{(1)}$, at the $k$-th iteration the algorithm first computes a Gauss-Newton iteration $\mbox{\boldmath $\hat{\vartheta}$}^{(k)}\in\mathbb{R}^{n_{\theta}}$. For every element of $\mbox{\boldmath $\hat{\vartheta}$}^{(k)}\in\mathbb{R}^{n_{\theta}}$ greater than a fix threshold $\epsilon_1>0$, the segment of $\mathcal{S}^{(k)}$ corresponding to that parameter is bisected: thus a new subdivision $\mathcal{S}^{(k+1)}$ is defined adding to $\mathcal{S}^{(k)}$ all the computed middle points. Finally are selected only those parameters which are greater than a fixed threshold $\epsilon_2>0$, and we indicate with $\Lambda^{(k)}$ this ensemble; the other parameters remain constant in the following iteration. The main ideas of the adaptive algorithm are sketched in algorithm \ref{algoritmo}.
\begin{algorithm}
	\footnotesize
	\caption{\small Sketch of the adaptive algorithm:}
\label{algoritmo}
  \begin{algorithmic}[1]
\STATE Given the finest subdivision of $\Gamma_h$ of step length $\Delta x$, the tolerance $tol>0$ and thresholds $\epsilon_1,\epsilon_2>0$, consider the coarse subdivision $\mathcal{S}^{(1)}=\left\{x^1_1,\ldots,x^1_{\frac{n^1_{\theta}}{2}+1}\right\}$, of $\Gamma_h$;
\STATE {$\mbox{\boldmath $\hat{\vartheta}$}^1=\textbf{0}_{n^1_{\theta}}\in\mathbb{R}^{n^1_{\theta}}$;}
\STATE{$l=1$, $\Lambda^{(1)}=[1,\ldots,n^1_{\theta}]$, set of indexes of parameters to be optimized}
\WHILE{$\mathcal{F}_d(\mbox{\boldmath $\hat{\vartheta}$}^{l}) < tol$}
\STATE{$\mathcal{S}^{(l+1)}:=\left\{x^{l+1}_1,\ldots,x^{l+1}_{\frac{n^{l+1}_{\theta}}{2}+1}\right\}=\mathcal{S}^{(l)}$;}
\STATE{$n^{l+1}_{\theta}=n^{l}_{\theta}$, $I=n^{l+1}_{\theta}$}
\FORALL{$i\in [1,I]$}
\IF{$\hat{\theta}^l(i) >\epsilon_1$\% \sffamily \scriptsize bisect the corresponding segment \normalfont\footnotesize}
\STATE{$n^{l+1}_{\theta}=n^{l+1}_{\theta}+1$, $I=I+1;$}
\STATE{let $[x^{l+1}(\hat{\theta}^l(i)),x^{l+1}(\hat{\theta}^l(i))]$ be the segment corresponding to parameter $\hat{\theta}^l(i)$;}
\STATE{$\mathcal{S}^{(l+1)}=\mathcal{S}^{(l+1)}\cup\frac{x^{l+1}(\hat{\theta}^l(i))-x^{l+1}(\hat{\theta}^l(i))}{2}$,}
\ENDIF
\ENDFOR
\STATE{given the subdivision $\mathcal{S}^{(l+1)}$ apply the projected damped Gauss-Newton method, optimizing \textbf{only} parameters whose indexes belong to $\Lambda^{(k)}$, obtaining $\mbox{\boldmath $\hat{\vartheta}$}^{l+1}\in\mathbb{R}^{n^{l+1}_{\theta}}$}
\STATE{$\Lambda^{(k)}=\emptyset$;}
\FORALL{$i\in [1,I]$}
\IF{$\hat{\vartheta}^{l+1}(i) >\epsilon_2$}
\STATE{$\Lambda^{(k)}=\Lambda^{(k)} \cup i$;}
\ENDIF
\ENDFOR
\STATE{$l=l+1$;}
\ENDWHILE
\end{algorithmic}		
\end{algorithm}

To avoid a large over-refinement, the bisection procedure can be limited, for example applying it at each iteration only a certain number of times, choosing the segments to be refined as those corresponding to greater parameters.

\subsubsection{Algorithm 4: using an adaptive parametrization and time localization}
The idea now is to reduce both the number of columns and of rows of the sensitivity matrix. As in algorithm 2, the domain $\Omega$ is partitioned in $n_s>1$ \textit{sections} $\mathcal{U}=\left\{s_j\right\}$, $j=1,\ldots,n_s$, however in algorithm 4 we assume that $\left\{\xi_1,\ldots,\xi_{n_s+1}\right\}=\mathcal{S}^{(1)}$, i.e. it coincides with the coarse initial subdivision applied in the adaptive strategy. In section $s_j$, considering the time interval $[t_{0}^{(j)},t_{f}^{(j)}]$, all parameters belonging to $\mathcal{O}^{(j)}\cup I^{(j)}$ will be estimated, and the adaptive procedure will be applied until a minimum is reached. Observe that this coincides with an internal loop: the ideas are summarized in algorithm \ref{time_loc_4}.
\begin{algorithm}
	\footnotesize
	\caption{\small Sketch of the adaptive algorithm with time localization:}
\label{time_loc_4}
  \begin{algorithmic}[1]
\STATE{Given the partition of $\Omega$ $\left\{\xi_1,\ldots,\xi_{n_s+1}\right\}=\mathcal{S}^{(1)}$, the thresholds $\epsilon_1,\epsilon_2,\epsilon_3>0$, $\mbox{\boldmath $\hat{\vartheta}$}^{0}=\textbf{0}$, $j=0$;}  
\WHILE{$\mathcal{F}_d(\mbox{\boldmath $\hat{\vartheta}$}^{j}) < tol$}  
\STATE{$j=j+1$;}
\STATE{$i=n_s$; $\mathcal{O}^{(n_s,1)}=\emptyset$}
\WHILE{$i>0$}
\STATE{$l=1$, $\Lambda^{(i,1)}=[1,\ldots,n^{j,i,1}_{\theta}]$, set of indexes of parameters to be optimized}
\WHILE{a minimum is reached}
\STATE {Let $I^{(i,l)}$ be the set of parameters of $\mbox{\boldmath $\hat{\vartheta}$}^{j,i,l}$ that belongs to section $i$;}
\STATE{in $[t_{0}^{(i)},t_{f}^{(i)}]$ apply the regularized projected damped Gauss Newton method to optimize parameters whose indices belong to $I^{(i,l)} \cup \mathcal{O}^{(i,l)}$;}
\STATE{update the positions $I^{(i,l)} \cup \mathcal{O}^{(i,l)}$ of $\mbox{\boldmath $\hat{\vartheta}$}^{j,i,l}$;}
\STATE{apply the adaptive strategy:}
\STATE{$\mathcal{S}^{(j,i,l+1)}=\mathcal{S}^{(j,i,l)}$;}
\STATE{$n^{j,i,l+1}_{\theta}=n^{j,i,l}_{\theta}$, $I=n^{j,i,l+1}_{\theta}$}
\FORALL{$k\in [1,I]$}
\IF{$\hat{\theta}^{j,i,l}(k) >\epsilon_1$}
\STATE{update $\mathcal{S}^{(j,i,l+1)}$, bisecting the segment corresponding to $\hat{\theta}^{j,i,l}(k)$;}
\ENDIF
\ENDFOR
\STATE{given the subdivision $\mathcal{S}^{(j,i,l+1)}$ apply the projected damped Gauss-Newton method, optimizing \textbf{only} parameters whose indexes belong to $\Lambda^{(j,i,l+1)}$}
\STATE{$\Lambda^{(j,i,l+1)}=\emptyset$;}
\FORALL{$k\in [1,I]$}
\IF{$\hat{\vartheta}^{j,i,l+1}(k) >\epsilon_2$}
\STATE{$\Lambda^{(j,i,l+1)}=\Lambda^{(j,i,l+1)} \cup k$;}
\ENDIF
\ENDFOR
\STATE{$l=l+1$;}
\STATE{if the subdivision has been refined, update $\mathcal{O}^{(i,l)}$}
\ENDWHILE
\STATE{$\mathcal{S}^{(j,i)}=\mathcal{S}^{(j,i,l)}$;}
\STATE{define $\mathcal{O}^{(j,i-1)}$ as the set of indices of parameters greater than $\epsilon_3$;}
\STATE{$i=i-1$;}
\ENDWHILE
\STATE{$\mbox{\boldmath $\hat{\vartheta}$}^{j}=\mbox{\boldmath $\hat{\vartheta}$}^{j,i,l}$}
\ENDWHILE
\end{algorithmic}		
\end{algorithm}

\subsubsection{Time localization: how to choose time intervals $[t_{0}^{(i)},t_{f}^{(i)}]$}
\begin{figure}[h]
\begin{center}
\includegraphics*[width=10cm]{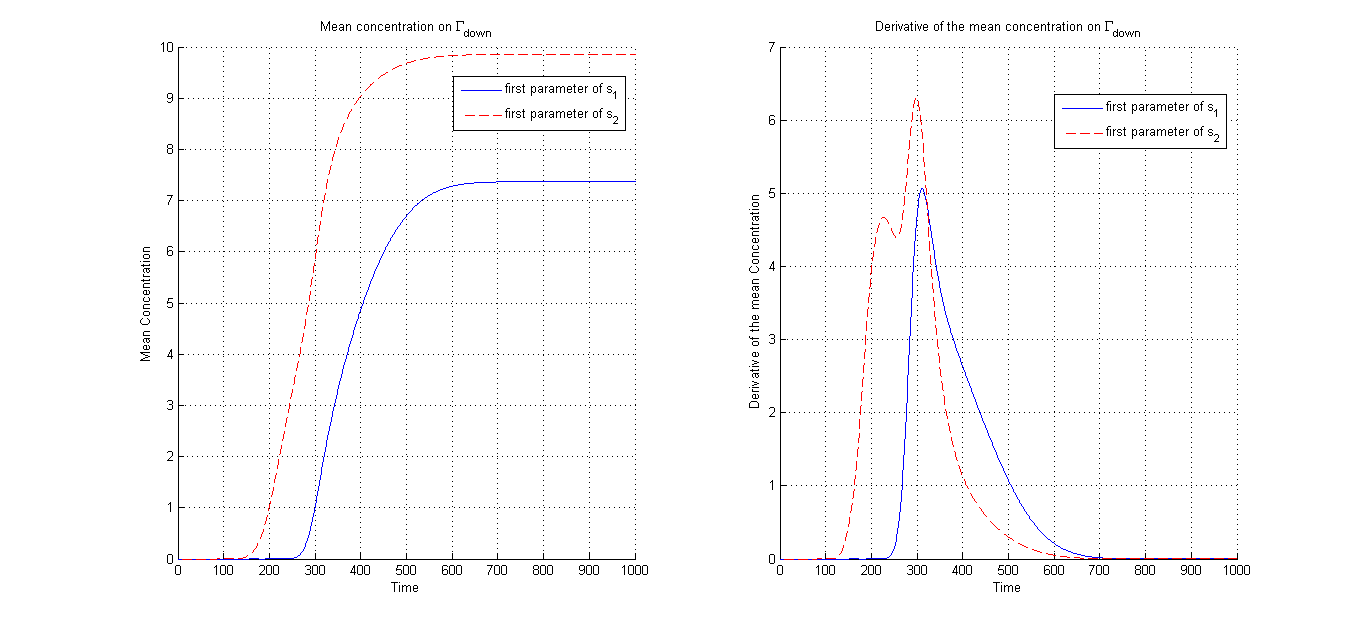}
\end{center}
\caption{\small \textit{Time evolution of the mean concentrations at $\Gamma_{down}$, $\zeta_1$ and $\zeta_2$, (left) and their derivative (right) for different boundary controls: the boundary control is different from zero only in the most left position of the finest subdivision of $s_1$ (blue) and $s_2$ (red).}}\normalsize
\label{mean_conc}
\end{figure}

A key point is the choice of the local time intervals $[t_{0}^{(i)},t_{f}^{(i)}]$, for every section $s_i$, $i=1,\ldots,n_s$, $t_{0}^{(i)}\geq t_0$ and $t_{f}^{(i)}\leq t_f$. 
\begin{figure}[h]
\begin{center}
\includegraphics*[width=10cm]{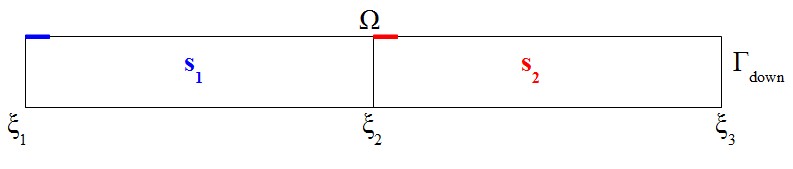}
\end{center}
\caption{\small \textit{Partition of $\Omega$ in $2$ sections to draw the curves of figure \ref{mean_conc}: to obtain the red (blue) curve of figure \ref{mean_conc}, it is considered the mean concentration on $\Gamma_{down}$, obtained imposing a control different from zero only in the most left segment of the finest subdivision of the upper horizontal segment of section $s_2$ ($s_1$), indicated in red (blue).}}\normalsize
\label{time_loc_dominio}
\end{figure}
The $i$-th interval must be chosen such that it contains the transitional dynamics of section $s_i$ but not that of sections $s_{j}$, $j<i$. To describe more clearly this idea, consider the model problem introduced in section \ref{numerical_results_known}: moreover suppose for simplicity that $n_s=2$, $\left\{\xi_1,\xi_{2},\xi_{3}\right\}=\left\{0,4,8\right\}$, as depicted in figure \ref{time_loc_dominio}, and consider as the finest subdivision a uniform one of step length $0.5$. Consider figure \ref{mean_conc}: the $j$-th curve $\zeta_j$, $j=1,2$, represents the mean concentration (left) and its derivative (right) at the outflow when the boundary control is different from zero only in the most left position of $s_j$ with respect to the finest subdivision. 
The interval corresponding to $s_2$ can be $[t_{0}^{(2)},t_{f}^{(2)}]=[180,260]$, when the red dotted curve corresponding to $s_2$, $\zeta_2$ , is increasing (transitional regime) and the blue curve corresponding to $s_1$, $\zeta_1$, is flat, i.e. when only the pollutant put into $\Omega$ in $s_2$ could reach $\Gamma_{down}$. While in $s_1$ the choice can be $[t_{0}^{(1)},t_{f}^{(1)}]=[240,400]$, since in this interval the transitional regime of $s_1$ occurs, as showed by $\zeta_1$. This intervals are used in section \ref{numerical_results}, to test algorithm 4.

The previous idea can be extended more rigorously to a general number of sections: let $\zeta_i$, $i=1,\ldots,n_s$, be the mean concentration at the outflow $\Gamma_{down}$ when the boundary control is different from zero only in the most left position of $s_i$, with respect to the finest subdivision. Consider a small threshold $\epsilon_4>0$, and two positive parameters $d,D >0$. Given 
\begin{equation*}
\begin{array}{r l l}
t_{0}^{(n_s)} & = & \min_{t\in [t_0,t_f]}\left\{\zeta_{n_s}^{'}(t)>\epsilon_4 \text{ and } \zeta_{n_s-1}^{'}(t)<\epsilon_4\right\},\\
t_{f}^{(n_s)} & = & \max_{t\in [t_0,t_f]}\left\{\zeta_{n_s}^{'}(t)>\epsilon_4 \text{ and } \zeta_{n_s-1}^{'}(t)<\epsilon_4\right\},
\end{array}
\end{equation*}
then for $i=1,\ldots,n_s-1$
\begin{equation*}
\begin{array}{r l l}
t_{0}^{(i)} & = & t_{f}^{(i+1)} - d,\\
t_{f}^{(i)} & = & \left\{\begin{array}{l l}\max_{t\in [t_0,t_f]}\left\{\zeta_i^{'}(t)>\epsilon_4 \text{ and } \zeta_{i-1}^{'}(t)<\epsilon_4\right\}, & \qquad \text{$i>1$}\\
\min\left\{t_f^{i+1}+D , \max_{t\in [t_0,t_f]}\left\{\zeta_i^{'}(t)>\epsilon_4 \text{ and } \zeta_{i-1}^{'}(t)<\epsilon_4\right\}\right\}, & \qquad \text{$i=1$}.\end{array}\right.
\end{array}
\end{equation*}

The parameter $d$ allows a small overlapping between local time intervals, while $D$ could limit the length of the inteval $[t_0^{(1)},t_f^{(1)}]$: e.g. in the example presented above, considering $n_s=2$, $d=0.2$ and $D=1.6$.


Observe that the definition of the intervals $[t_{0}^{(i)},t_{f}^{(i)}]$ depends on the shape of the domain, on the velocity field and on the coefficients of the PDE (\ref{direct_problem}): each time one of them is changed, also the intervals should be estimated, observing the transitional dynamics of each section, as explained above.

\subsection{Comparing computational costs}
In this section we compare the computational costs of the four algorithms.

The first one consists in using the finest subdivision, with the projected damped Gauss Newton strategy. The computational cost of each iteration is pretty high: since the solution of the direct problem has cost $N N_h^3$, computing the sensitivity matrix $\Psi_{\mbox{\boldmath $\vartheta$}}\in\mathbb{R}^{n_y N\times n^{(f)}_{\theta}}$ has cost $n^{(f)}_{\theta}N N_h^3$, where $n^{(f)}_{\theta}$ is the number of parameter of the finest subdivision, which is maximal. Moreover computing the SVD to obtain the new iteration has cost $4 n_y^2 N^2 n^{(f)}_{\theta} + 8 N n_y (n^{(f)}_{\theta})^2+9 (n^{(f)}_{\theta})^3$. Finally computing the new prediction error has cost $N N_h^3$.

To decrease the cost, the idea is to consider a sensitivity matrix of lower dimensions. The second algorithm consists in combining the finest subdivision with localization in time. The number of sections in this case coincides with one half of the number of parameters of the finest subdivision $n^{(f)}_{\theta}$. At each iteration $k$, for every section $i=1,\ldots,n_s$, $n_s=\frac{n^{(f)}_{\theta}}{2}$, computing $\Psi^{(i)}_{\mbox{\boldmath $\vartheta$}}\in\mathbb{R}^{n_y \frac{t_f^{(i)}-t_0^{(i)}}{Dt}\times n^{(k,i)}_{\theta}}$ costs $n^{(k,i)}_{\theta} (\frac{t_f^{(i)}-t_0^{(i)}}{Dt}) N_h^3$, where $n^{(k,i)}_{\theta}$ denotes the cardinality of $I^{(i)}\cup\mathcal{O}^{(i)}$. Moreover computing the SVD to obtain the new iteration has cost $4 n_y^2 \left(\frac{t_f^{(i)}-t_0^{(i)}}{Dt}\right)^2 n^{(k,i)}_{\theta} + 8 \frac{t_f^{(i)}-t_0^{(i)}}{Dt} n_y (n^{(k,i)}_{\theta})^2+ (n^{(k,i)}_{\theta})^3$. Finally computing the new prediction error has cost $N N_h^3$. Although an higher number of systems must be solved, the algorithm is less costly since the sensitivity matrix has much lower dimensions.

Another possibility to decrease the cost of algorithm one, is to use the third algorithm, which consists in adopting an adaptive parametrization. At the $k-$th iteration computing the sensitivity matrix $\Psi_{\mbox{\boldmath $\vartheta$}}\in\mathbb{R}^{n_y N\times n^{(k)}_{\theta}}$ has cost $n^{(k)}_{\theta}N N_h^3$, where the number of parameter $n^{(k)}_{\theta}$ varies during the iterations and $n^{(k)}_{\theta}<n^{(f)}_{\theta}$. Moreover computing the SVD to obtain the new iteration has cost $4 n_y^2 N^2 n^{(k)}_{\theta} + 8 N n_y (n^{(k)}_{\theta})^2+9 (n^{(k)}_{\theta})^3$. Finally computing the new prediction error has cost $N N_h^3$. The gain with respect to the first strategy is evident if $n^{(k)}_{\theta}<<n^{(f)}_{\theta}$. 

The fourth algorithm combines both time localization and the adaptive parametrization. The number of sections in this case coincides with one half the number of parameters of the initial coarse subdivision $\mathcal{S}^{(1)}$. The difference with respect to the second algorithm is that the number of sections $n_s$ is lower, because it is no more related to the finest subdivision: in fact the adaptive parametrization guides the choice of parameters to be estimated at each iteration. However the introduction of the adaptive parametrization introduces an inner loop. In detail, at each iteration $k$, for every section $i=1,\ldots,n_s$, applying the adaptive procedure until a minimum is reached (index $l$), computing $\Psi^{(k,i,l)}_{\mbox{\boldmath $\vartheta$}}\in\mathbb{R}^{n_y \frac{t_f^{(i)}-t_0^{(i)}}{Dt}\times n^{(k,i,l)}_{\theta}}$ costs $n^{(k,i,l)}_{\theta} (\frac{t_f^{(i)}-t_0^{(i)}}{Dt}) N_h^3$. Moreover computing the SVD to obtain the new iteration has cost $4 n_y^2 \left(\frac{t_f^{(i)}-t_0^{(i)}}{Dt}\right)^2 n^{(k,i,l)}_{\theta} + 8 \frac{t_f^{(i)}-t_0^{(i)}}{Dt} n_y (n^{(k,i,l)}_{\theta})^2+ (n^{(k,i,l)}_{\theta})^3$. Finally computing the new prediction error has cost $N N_h^3$.

Just to give an idea of the computational gain of the fourth algorithm, computational costs of the four algorithms are summarized in table \ref{table_results_fiume_costo}, averaging results of tests presented in section \ref{numerical_results}.
\begin{center}
\begin{table}[h]
\tiny{\begin{tabular}{|c||c|c|c|c|} \hline
\multicolumn{1}{|c||}{} & \multicolumn{1}{|c|}{\textbf{Finest subdivision}} & \multicolumn{1}{|c|}{ \textbf{Finest subdivision}} & \multicolumn{1}{|c|}{\textbf{Adaptive subdivision}} &  \multicolumn{1}{|c|}{ \textbf{Adaptive subdivision}} \\
 \multicolumn{1}{|c||}{} & \multicolumn{1}{|c|}{} & \multicolumn{1}{|c|}{\textbf{+ time localization}} & \multicolumn{1}{|c|}{} &  \multicolumn{1}{|c|}{\textbf{+ time localization}} \\ \hline
\textbf{Computational} &  &  &  & \\ 
\textbf{cost} & $8 \cdot 10^{14}$ & $2 \cdot 10^{13}$ & $5 \cdot 10^{12}$ & $8 \cdot 10^{11}$\\ \hline
\end{tabular}}
\caption{\small \textit{Estimated computational cost of the four algorithms: using the finest subdivision and the projected damped Gauss Newton method, using the finest subdivision and the localization in time, using the adaptive parametrization and using the adaptive parametrization and time localization.}}  
\label{table_results_fiume_costo}
\end{table}
\end{center}

\subsection{Numerical results}
\label{numerical_results}
\begin{figure}[t]
\begin{center}
\includegraphics*[width=6cm]{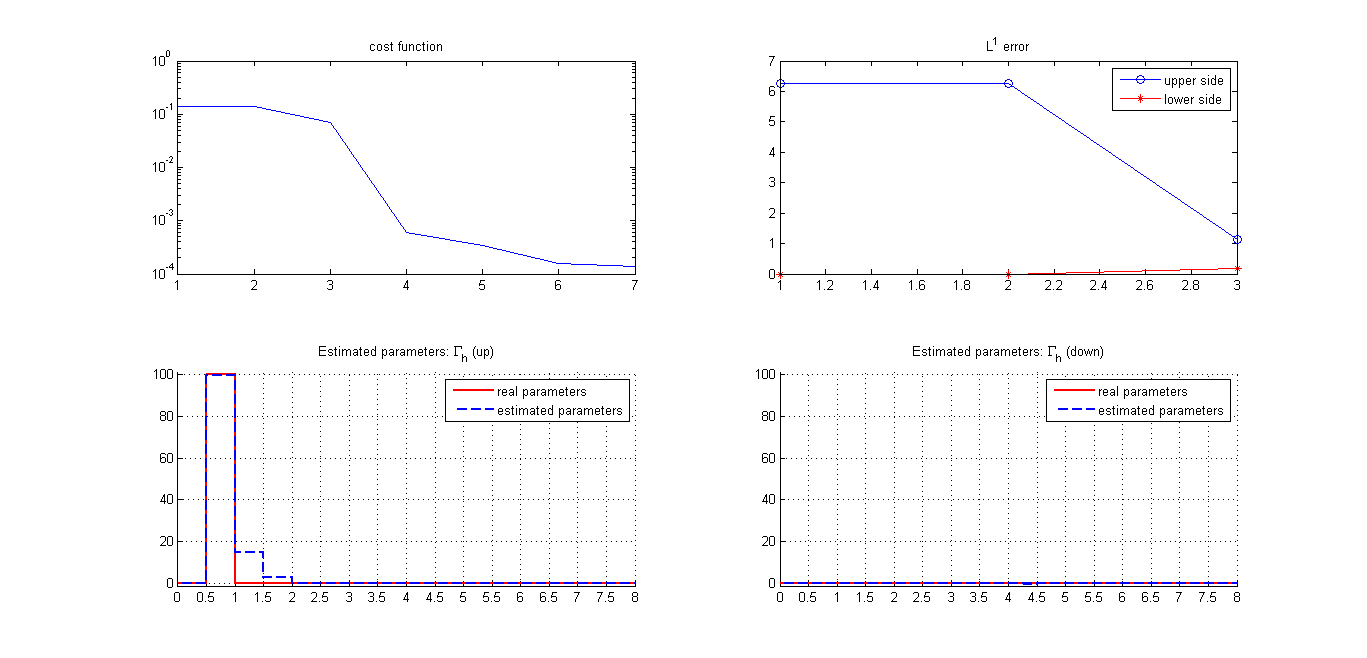}
\includegraphics*[width=6cm]{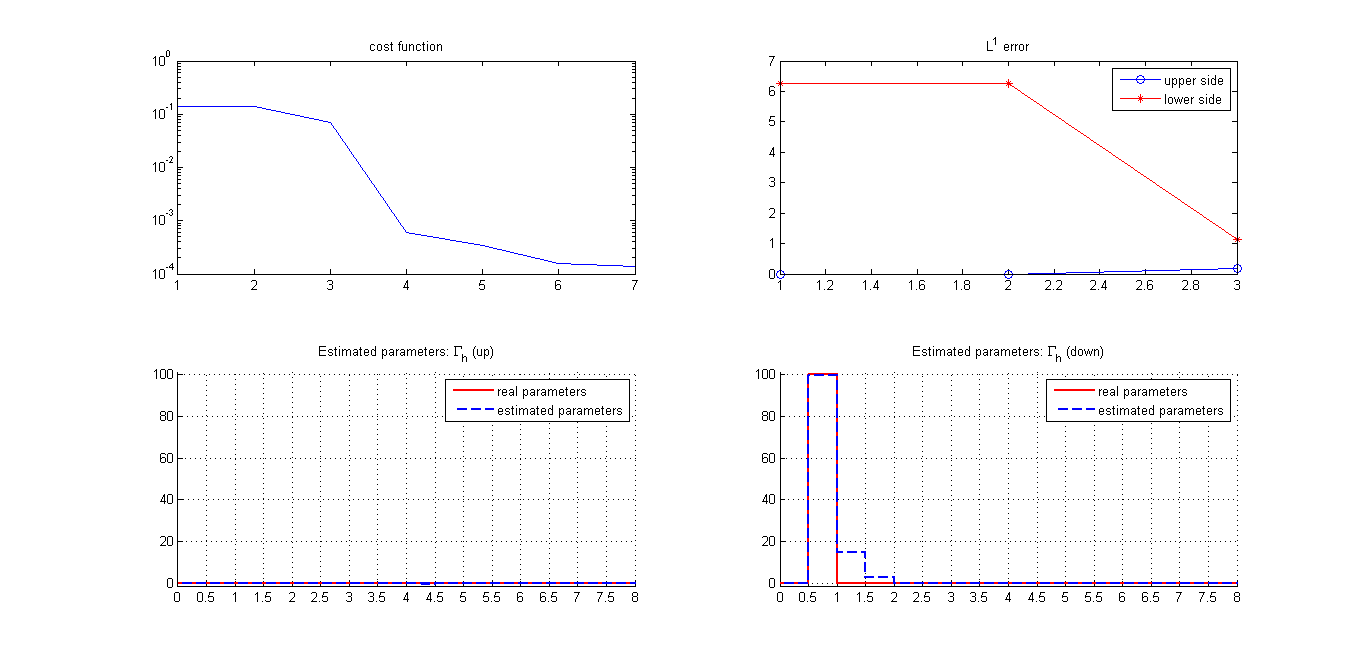}
\includegraphics*[width=6cm]{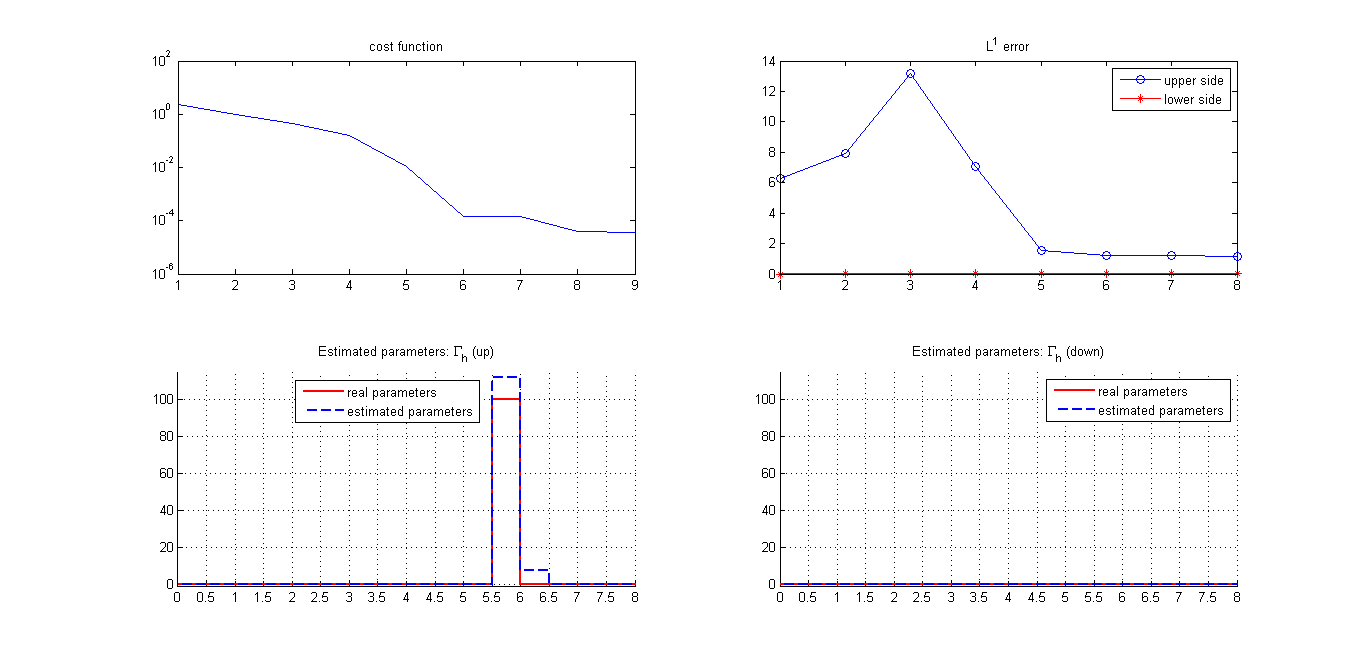}
\includegraphics*[width=6cm]{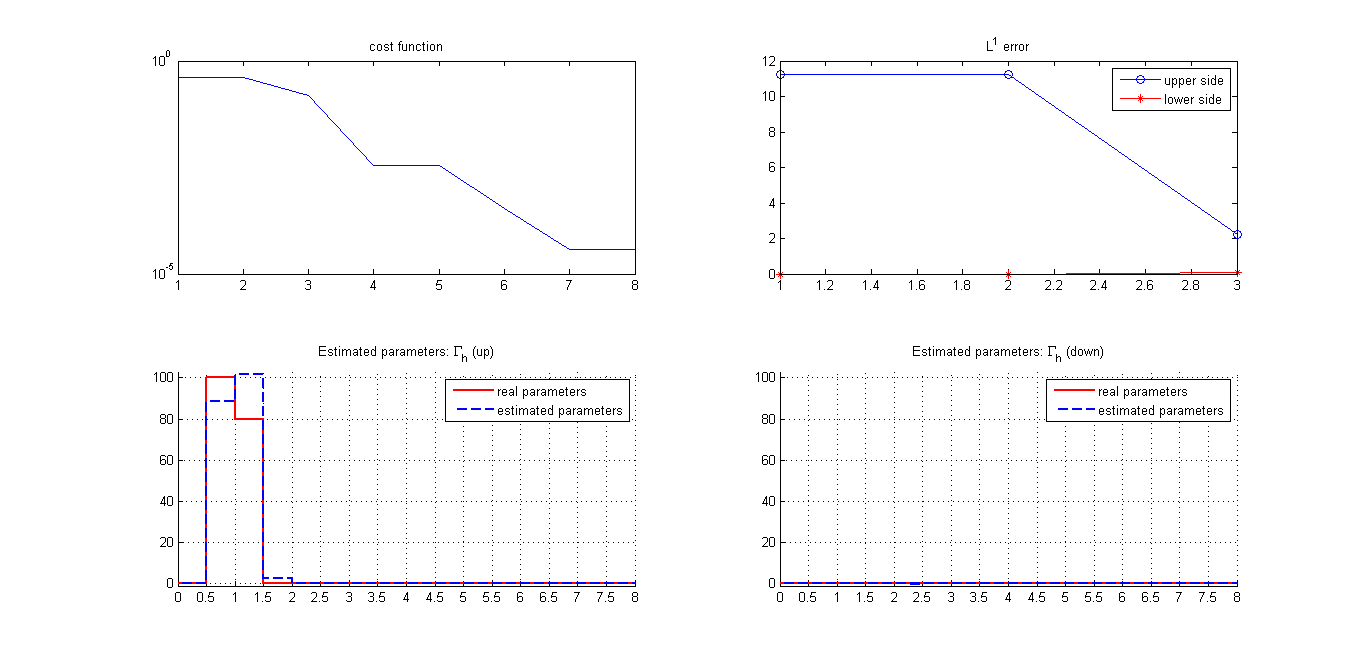}
\includegraphics*[width=6cm]{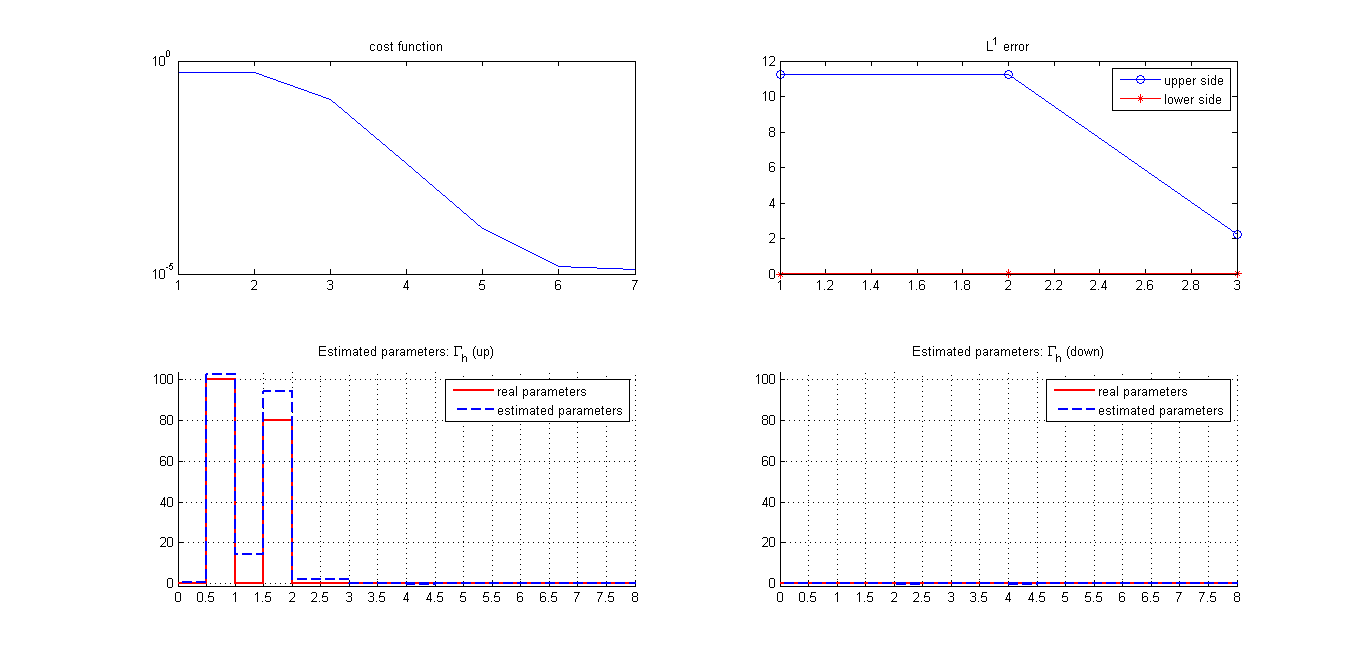}
\includegraphics*[width=6cm]{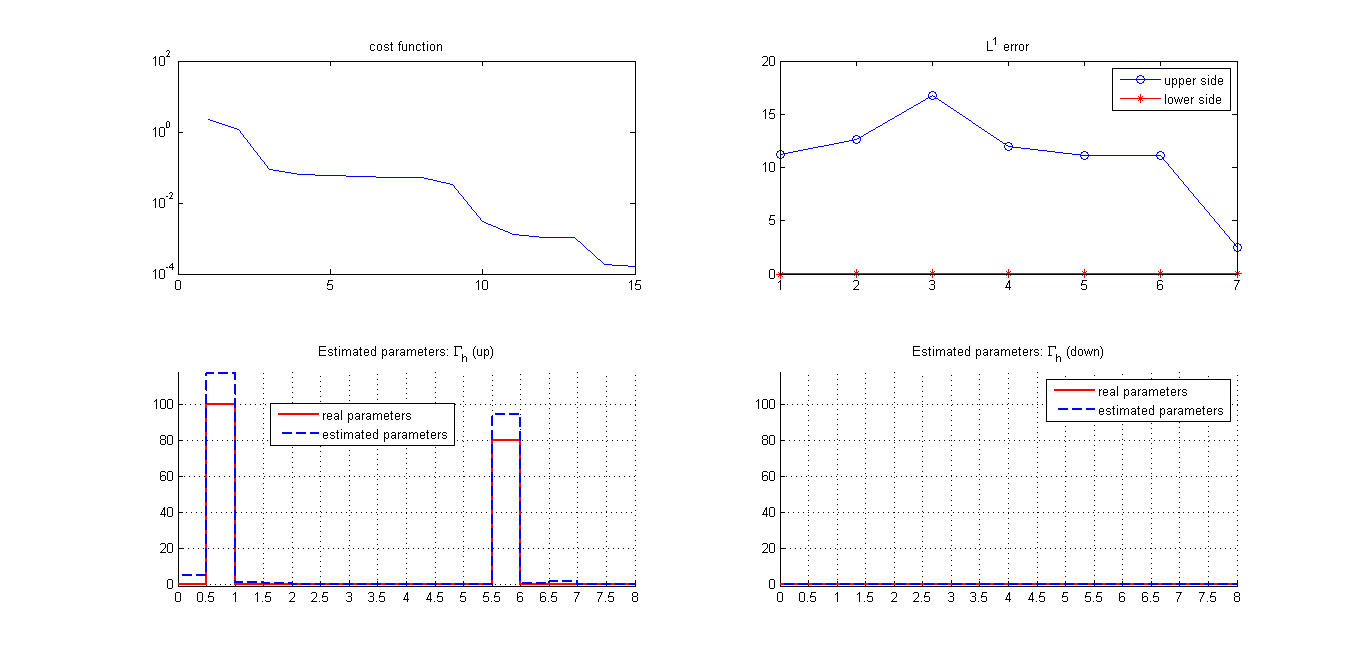}
\includegraphics*[width=6cm]{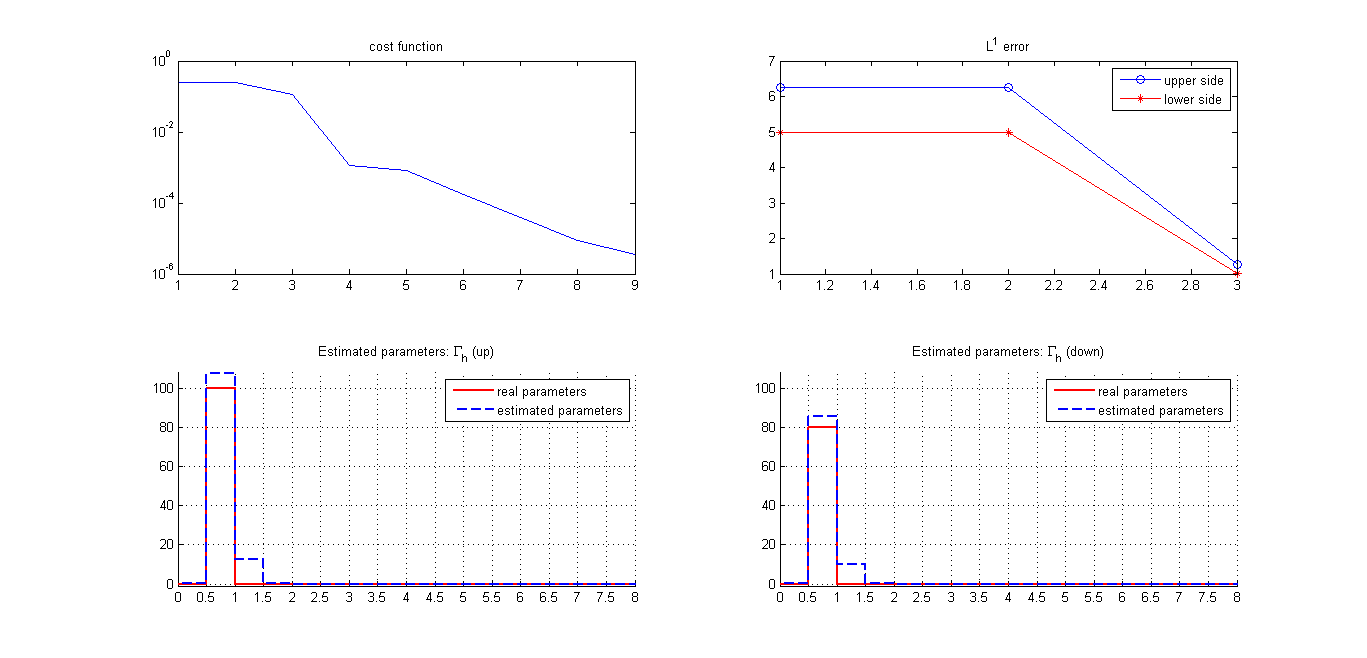}
\includegraphics*[width=6cm]{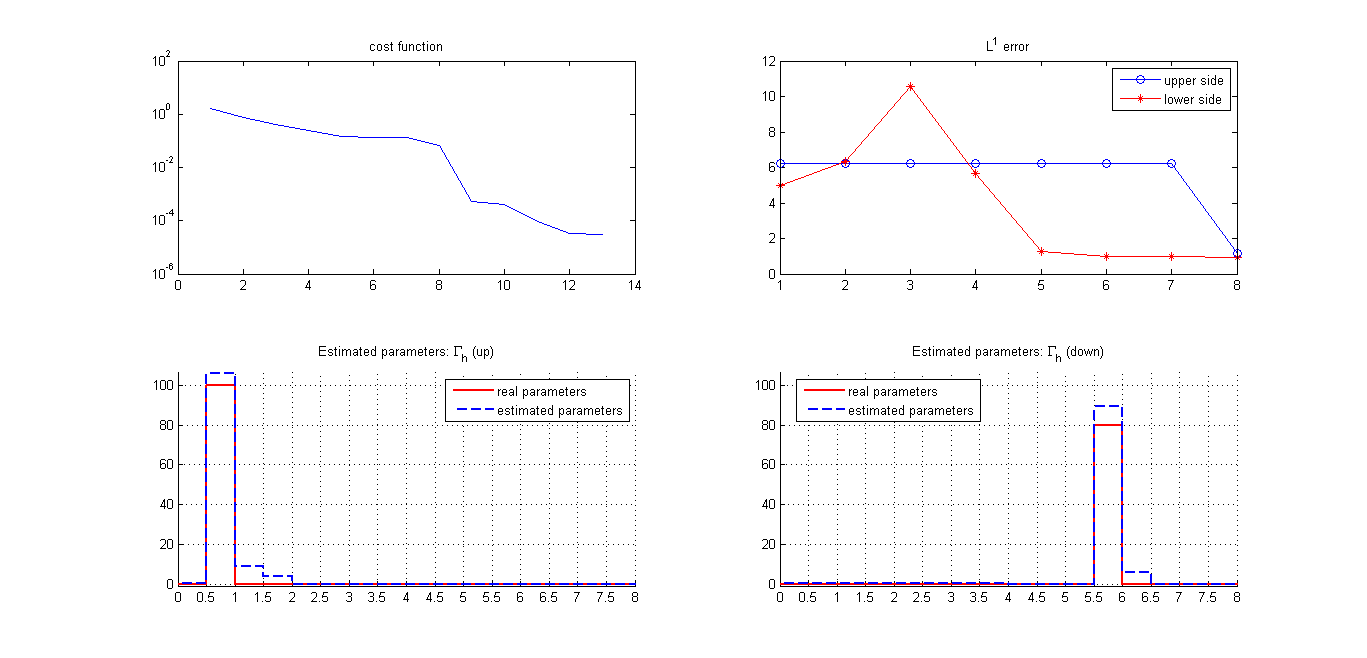}
\includegraphics*[width=6cm]{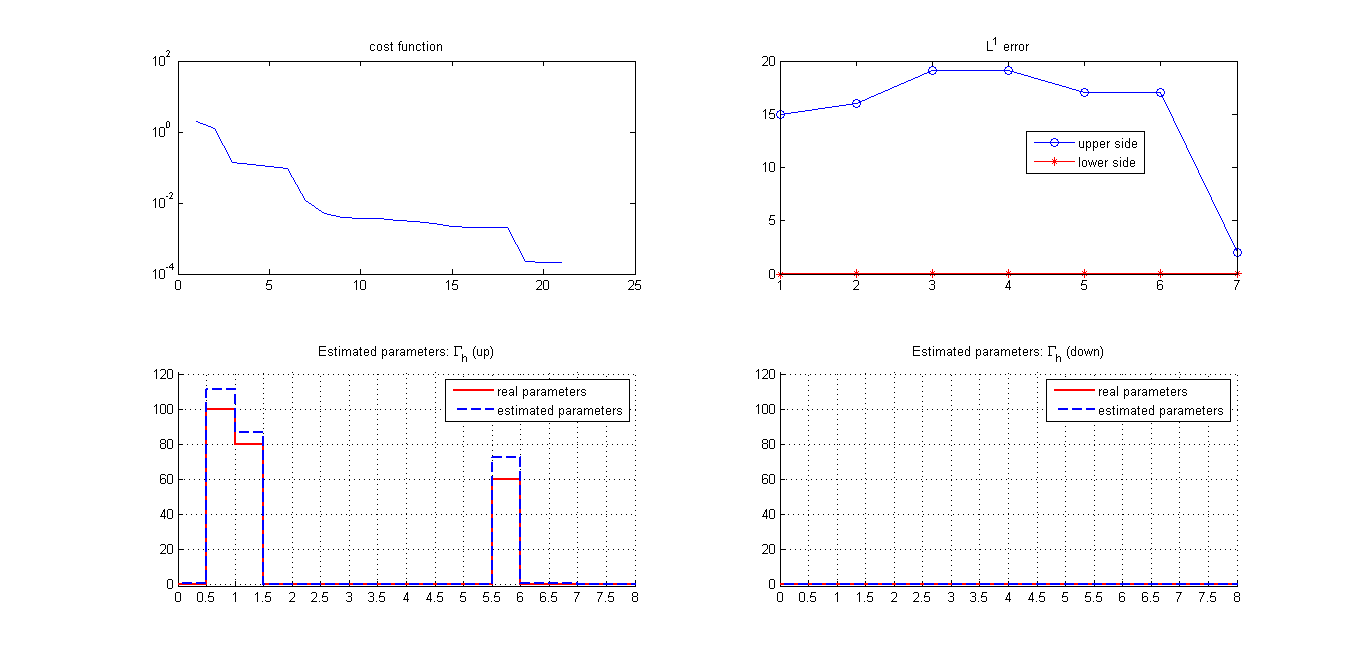}
\end{center}
\caption{\small \textit{Nine test cases: results of the adaptive strategy with time localization: computed estimate (blu dotted line), real control (red line). For each figure: cost function (first row, left), $L^1$ error (first row, right), approximation of the upper horizontal segment (second row, left), approximation of the bottom horizontal segment (second row, right).}}\normalsize
\label{num_res_fiume}
\end{figure}

In this section we present some numerical tests to verify the effectiveness of the algorithm. As in section \ref{numerical_results_known}, experimental data are simulated numerically, on $\Omega = [0,8]\times [0,1]$,  $\Gamma_h=[0,8]\times\left\{1\right\}\cup [0,8]\times\left\{0\right\}$. Moreover the velocity field $\textbf{u}$ is modeled as a Poiseuille flow i.e.
$$\textbf{u}(x_1,x_2)=\left(\begin{array}{c c} -4\nu x^2_2 + 4 \nu x_2\\ 0 \end{array}\right).$$
We assume that $\mu=0.1$, $\sigma=0.1$ and $C_{up}=0.1$. Moreover we consider the finest subdivision with step length $\Delta x =0.5$. In algorithm 2 we consider $\left\{\xi_1,\ldots,\xi_{n_s+1}\right\}$ coincident with the finest subdivision, while in algorithm 4 $n_s=2$ and $\left\{\xi_1,\ldots,\xi_{n_s+1}\right\}=\left\{0,4,8\right\}$. Define \textit{optimal subdivision} the one which describes the real profile with the minimum number of parameters using the bisection criterium. With \textit{distance from the optimal subdivision} we indicate the number of points added (sign +) or subtracted (sign -) to the optimal sundivision.
We consider 9 test cases: results for the adaptive strategy with localization in time are shown in figure \ref{num_res_fiume}. In table \ref{table_results_fiume}, four algorithms are compared: using the finest subdivision and the projected damped Gauss Newton method, using the finest subdivision and the localization in time, using the adaptive parametrization and using the adaptive parametrization and time localization.

First of all observe that the number of iterations of algorithms 2 and 4 is higher since also sub-iterations to reach the minimum inside each section are counted (inner loop).

In tests 1, 2 and 7, also working on the finest subdivision performs well, but it is much more costly. When the condition number of the sensitivity matrix $\Psi_{\mbox{\boldmath $\vartheta$}}$ increases, the accuracy is low. In particular in tests 3, 4, 5 it is evident how time localization improves convergence results both in algorithms 2 and 4, selecting only some rows of $\Psi_{\mbox{\boldmath $\vartheta$}}$. However adopting only time localization is not sufficient in tests 6,7,9. Using an adaptive parametrization corresponds to select only some columns of $\Psi_{\mbox{\boldmath $\vartheta$}}$, considering a less number of parameters: in algorithm 3 the number of points added to the optimal subdivision is very low, but in general the estimates tend to be too much approximated. The best strategy consists in combining both adaptive parametrization and time localization (algorithm 4): this is a good compromise between good estimates and reasonable computational cost. Its effectiveness is evident e.g. in tests 8 and 9. Moreover it only adds few points to the optimal subdivision. 

\begin{center}
\begin{table}[t]
\scriptsize{\begin{tabular}{|c|c||c|c|c||c|c|c||c|c|c||c|c|c|} \hline
\multicolumn{2}{|c||}{\textbf{Test}} & \multicolumn{3}{|c||}{\textbf{Finest subdivision}} & \multicolumn{3}{|c||}{ \textbf{Finest subdivision}} & \multicolumn{3}{|c||}{\textbf{Adaptive subdivision}} &  \multicolumn{3}{|c|}{ \textbf{Adaptive subdivision}} \\
 \multicolumn{2}{|c||}{} & \multicolumn{3}{|c||}{} & \multicolumn{3}{|c||}{\textbf{+ time localization}} & \multicolumn{3}{|c||}{} &  \multicolumn{3}{|c|}{\textbf{+ time localization}} \\ \hline
 &  & \textit{up} &  \textit{down} & & \textit{up} & \textit{down} & & \textit{up} &  \textit{down} & & \textit{up} &  \textit{down} & \\ \hline\hline

\textbf{1} & \textit{$L^1$-err} & $10^{-12}$ & $10^{-12}$ & & 0.442 & 0.442 & & 7.69 & 0.12 & & 1.15 & 0.168 & \\ \hline
					 & \textit{opt. sub.} & +11 & +14 & & +11 & +14 & & +1 & 0 & & +1 & 0 & \\ \hline
					 & \textit{num. it.} & & & 2 & & & 18 & & & 4 & & & 7 \\ \hline 
					 & $\mathcal{F}_d(\mbox{\boldmath $\vartheta$})$ & & & $10^{-20}$ & & & $10^{-6}$ & & & $10^{-5}$ & & & $10^{-4}$\\ \hline\hline

\textbf{2} & \textit{$L^1$-err} & $10^{-12}$ & $10^{-12}$ & & 0.02 & 0.02 & & 0.12 & 7.69 & & 0.168 & 1.15 & \\ \hline
					 & \textit{opt. sub.} & +14 & +11 & & +14 & +11 & & 0 & +1 & & 0 & +1 & \\ \hline
					 & \textit{num. it.} & & & 2 & & & 20 & & & 4 & & & 7 \\ \hline 
					 & $\mathcal{F}_d(\mbox{\boldmath $\vartheta$})$ & & & $10^{-20}$ & & & $10^{-6}$ & & & $10^{-5}$ & & & $10^{-4}$\\ \hline\hline

\textbf{3} & \textit{$L^1$-err} & 2.72 & 0.3974 & & 0 & 0 & & 1.33 & 0.02 & & 0.18 & $10^{-3}$ & \\ \hline
					 & \textit{opt. sub.} & +11 & +14 & & +11 & +14 & & 0 & +1 & & +2 & 0 & \\ \hline
					 & \textit{num. it.} & & & 6 & & & 17 & & & 9 & & & 9 \\ \hline 
					 & $\mathcal{F}_d(\mbox{\boldmath $\vartheta$})$ & & & $0.028$ & & & $10^{-20}$ & & & $0.0012$ & & & $10^{-5}$\\ \hline\hline

\textbf{4} & \textit{$L^1$-err} & 8.911 & 0.1102 & & 1.021 & $10^{-3}$ & & 11.25 & 0.16 & & 2.247 & 0.07 & \\ \hline
					 & \textit{opt. sub.} & +10 & +14 & & +10 & +14 & & 0 & 0 & & 0 & +1 & \\ \hline
					 & \textit{num. it.} & & & 5 & & & 13 & & & 4 & & & 8 \\ \hline 
					 & $\mathcal{F}_d(\mbox{\boldmath $\vartheta$})$ & & & $10^{-3}$ & & & $10^{-3}$ & & & $10^{-3}$ & & & $10^{-5}$\\ \hline\hline
					 
\textbf{5} & \textit{$L^1$-err} & 5.576 & 0.047 & & 1.611 & $10^{-4}$ & & 12 & 0.14 & & 2.224 & 0.03 & \\ \hline
					 & \textit{opt. sub.} & +10 & +14 & & +10 & +14 & & -1 & 0 & & +1 & +1 & \\ \hline
					 & \textit{num. it.} & & & 5 & & & 11 & & & 3 & & & 7 \\ \hline 
					 & $\mathcal{F}_d(\mbox{\boldmath $\vartheta$})$ & & & $10^{-3}$ & & & $10^{-4}$ & & & $10^{-4}$ & & & $10^{-5}$\\ \hline\hline					 
					 
\textbf{6} & \textit{$L^1$-err} & 2.653 & 0.2871 & & 8.591 & $10^{-13}$ & & 2.33 & 0.01 & & 2.48 & 0.01 & \\ \hline
					 & \textit{opt. sub.} & +8 & +14 & & +8 & +14 & & +1 & 0 & & +3 & 0 & \\ \hline
					 & \textit{num. it.} & & & 5 & & & 17 & & & 10 & & & 15 \\ \hline 
					 & $\mathcal{F}_d(\mbox{\boldmath $\vartheta$})$ & & & $10^{-2}$ & & & 0.068 & & & $10^{-4}$ & & & $10^{-4}$\\ \hline\hline					 

\textbf{7} & \textit{$L^1$-err} & $10^{-13}$ & $10^{-13}$ & & 0.36 & 0.36 & & 7.63 & 6.12 & & 1.267 & 1.019 & \\ \hline
					 & \textit{opt. sub.} & +9 & +9 & & +9 & +9 & & 0 & 0 & & +1 & +1 & \\ \hline
					 & \textit{num. it.} & & & 2 & & & 17 & & & 4 & & & 9 \\ \hline 
					 & $\mathcal{F}_d(\mbox{\boldmath $\vartheta$})$ & & & $10^{-20}$ & & & $10^{-6}$ & & & $10^{-5}$ & & & $10^{-6}$\\ \hline\hline					 

\textbf{8} & \textit{$L^1$-err} & 1.969 & 1.002 & & 6.25 & 9.17 & & 14.9 & 8.9 & & 0.95 & 0.95 & \\ \hline
					 & \textit{opt. sub.} & +9 & +9 & & +9 & +9 & & +2 & 0 & & +1 & +2 & \\ \hline
					 & \textit{num. it.} & & & 5 & & & 17 & & & 5 & & & 13 \\ \hline 
					 & $\mathcal{F}_d(\mbox{\boldmath $\vartheta$})$ & & & $10^{-2}$ & & & 0.13 & & & 0.19 & & & $10^{-5}$\\ \hline\hline
			
\textbf{9} & \textit{$L^1$-err} & 14.22 & 0.1818 & & 14.34 & $10^{-12}$ & & 2.65 & 0.9 & & 2.01 & 0.004 & \\ \hline
					 & \textit{opt. sub.} & +7 & +14 & & +7 & +14 & & +3 & 0 & & +2 & 0 & \\ \hline
					 & \textit{num. it.} & & & 11 & & & 19 & & & 16 & & & 21 \\ \hline 
					 & $\mathcal{F}_d(\mbox{\boldmath $\vartheta$})$ & & & $10^{-2}$ & & & 0.2 & & & 0.001 & & & $10^{-4}$\\ \hline
					 
\end{tabular}}
\caption{\small \textit{Comparison between four algorithms: using the finest subdivision and the projected damped Gauss Newton method, using the finest subdivision and the localization in time, using the adaptive parametrization and using the adaptive parametrization and time localization. $L^1$-error in the upper and lower horizontal segments, number of points added to the optimal subdivision in the upper and lower horizontal segments, number of iterations and final cost function.}}  
\label{table_results_fiume}
\end{table}
\end{center}

\begin{figure}[h]
\begin{center}
\includegraphics*[width=3.5cm]{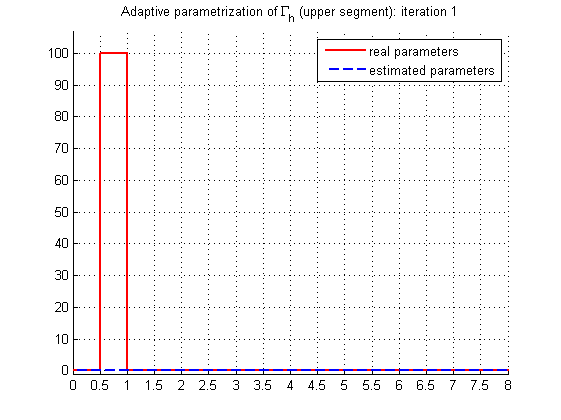}
\includegraphics*[width=3.5cm]{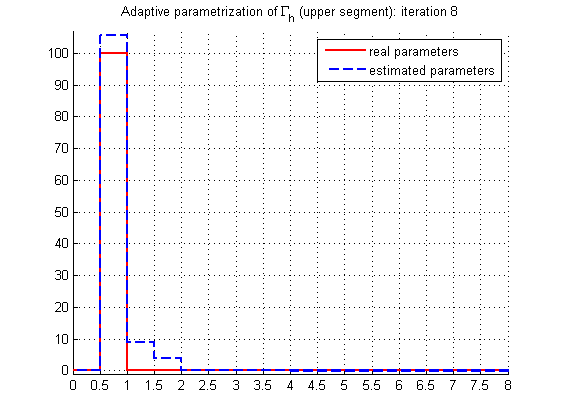}
\end{center}
\caption{\small \textit{Test 8. Adaptive parametrization and time localization. Evolution of the approximation (blu dotted line), real control (red line). Upper horizontal segment.}}\normalsize
\label{adaptive10_up}
\end{figure}

\begin{figure}[h]
\begin{center}
\includegraphics*[width=3.5cm]{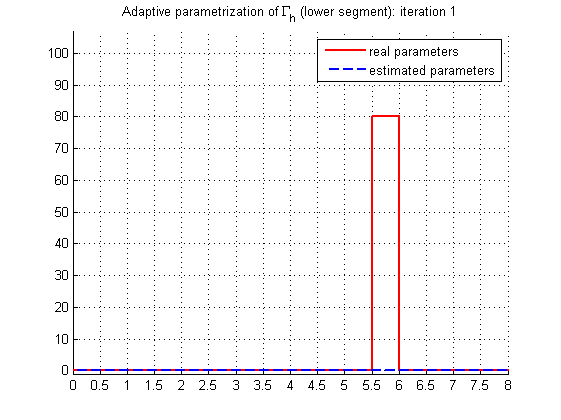}
\includegraphics*[width=3.5cm]{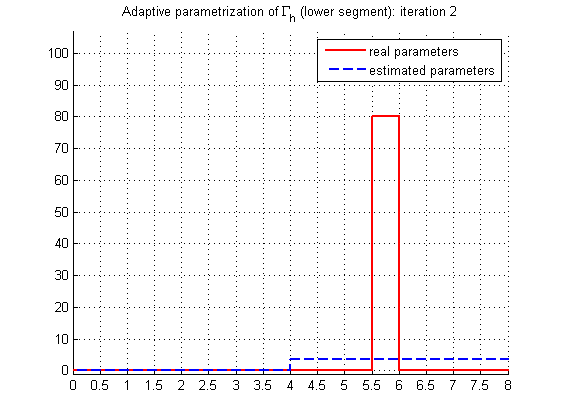}
\includegraphics*[width=3.5cm]{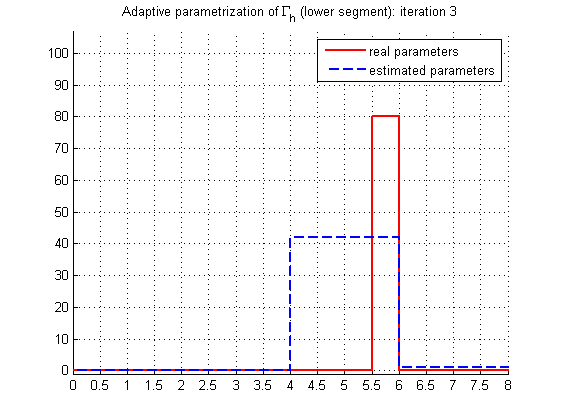}
\includegraphics*[width=3.5cm]{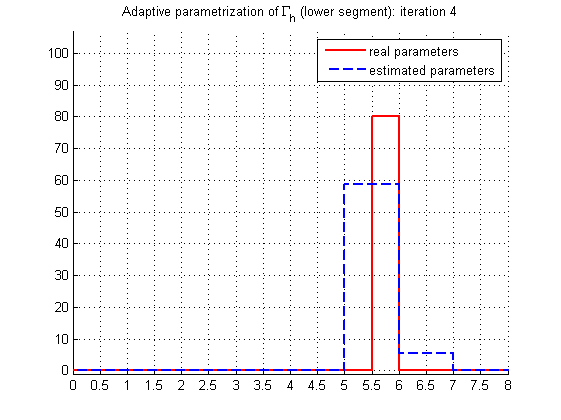}
\includegraphics*[width=3.5cm]{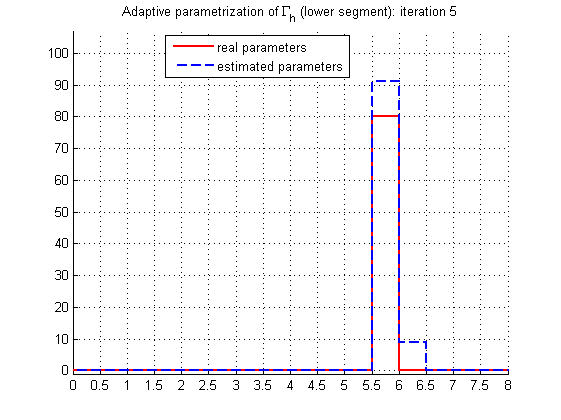}
\includegraphics*[width=3.5cm]{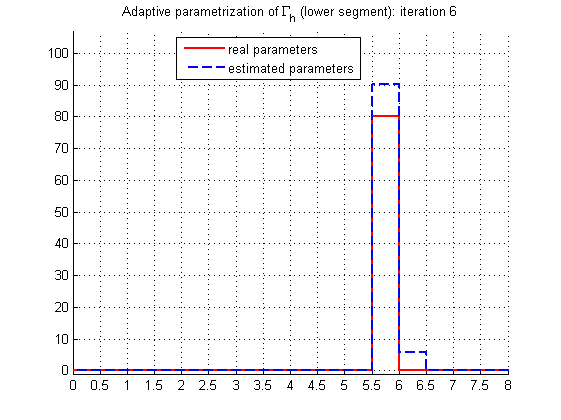}
\includegraphics*[width=3.5cm]{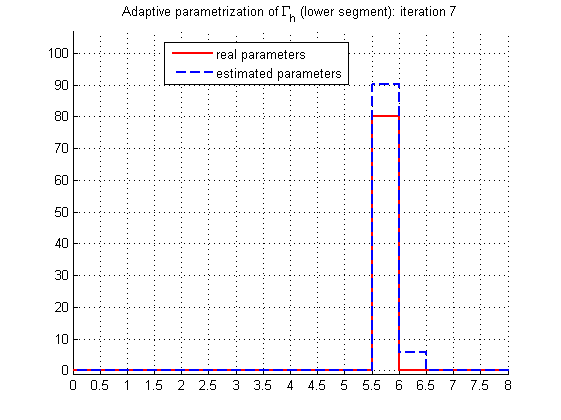}
\includegraphics*[width=3.5cm]{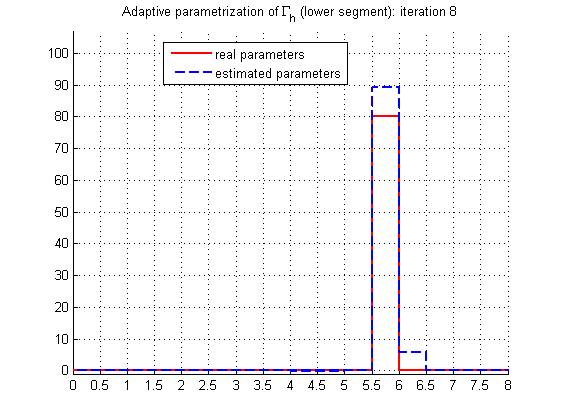}
\end{center}
\caption{\small \textit{Test 8. Adaptive parametrization and time localization. Evolution of the approximation (blu dotted line), real control (red line). Bottom horizontal segment.}}\normalsize
\label{adaptive10_down}
\end{figure}

\begin{figure}[h]
\begin{center}
\includegraphics*[width=10cm]{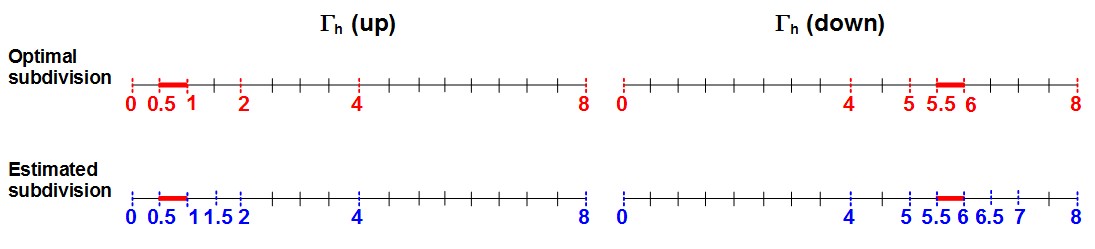}
\end{center}
\caption{\small \textit{Test 8. First row: optimal subdivision that could be obtained using a bisection strategy. Second row: estimated subdivision.}}\normalsize
\label{suddivisione_10}
\end{figure}

In figure \ref{adaptive10_up} and \ref{adaptive10_down} different iterations of the algorithm are shown for test 8: it is evident how the algorithm firstly optimize parameters of section $s_2=[4,8]\times [0,1]$ (figure \ref{adaptive10_down}), and then that of $s_1=[0,4]\times [0,1]$ (in figure \ref{adaptive10_up} the first 7 iteration are identical to the first one, thus only iterations 1 and 8 are plotted). The estimated subdivision is sketched in figure \ref{suddivisione_10}: it is evident how algorithm 4 slightly over-refine the optimal subdivision.


\subsection{Conditioning of the problem}
\begin{figure}[h]
\begin{center}
\includegraphics*[width=4cm]{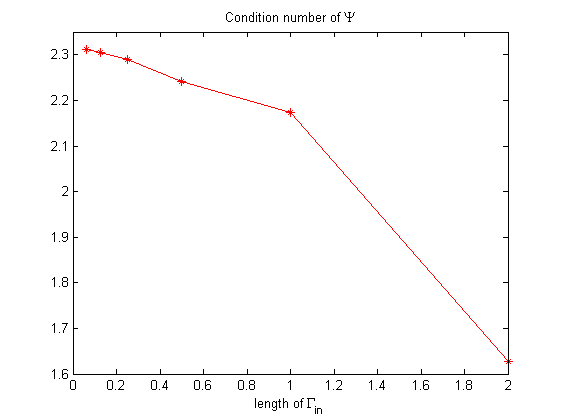}
\end{center}
\caption{\small \textit{Example 2, with $\Gamma_{in}=[5-h,5]\times\left\{1\right\}\cup[2-h,2]\times\left\{0\right\}$. Condition number of $\Psi_{\mbox{\boldmath $\hat{\vartheta}$}}$ for different values of $h=\left\{0.0625,0.125,0.25,0.5,1,2\right\}$.}}\normalsize
\label{cond_vs_length}
\end{figure}

The ill-conditioning of the system matrix $\Psi_{\mbox{\boldmath $\hat{\vartheta}$}}$ could increase when smaller segments are considered in $\Gamma_h$: in fact in this case consecutive columns tend to be close to linear dependence, due to the small distance ($\Delta x$) of the corresponding nodes in $\Gamma_h$. This can be demonstrated numerically: consider in fact example 2 presented in section \ref{numerical_results_known} and generalize it considering the following parametric problem
$$\Gamma_{in}=[5-h,5]\times\left\{1\right\}\cup[2-h,2]\times\left\{0\right\},\ \mbox{\boldmath $\vartheta$} = (100,80),\qquad 0 < h \leq 2.$$
Even supposing to know source localition $\Gamma_{in}$, solving the problem for different values of $h=\left\{0.0625,0.125,0.25,0.5,1,2\right\}$ and computing the condition number of the sensitivity matrix, it can be seen that as $h$ decreases, the condition number increases (cfr. figure \ref{cond_vs_length}).
Since the condition number of the sensitivity matrix could become higher when smaller segments are considered, working on the finest subdivision could not be effective to reduce the ill-conditioning of the problem and an adaptive parametrization should be preferred. Observe moreover that in adaptive algorithms the Gauss Newton method is applied only to those parameters belonging to $\Lambda^{(k)}$: avoiding parameters less than the threshold $\epsilon_2$ is useful to reduce columns linear dependence.

\begin{figure}[h]
\begin{center}
\includegraphics*[width=5cm]{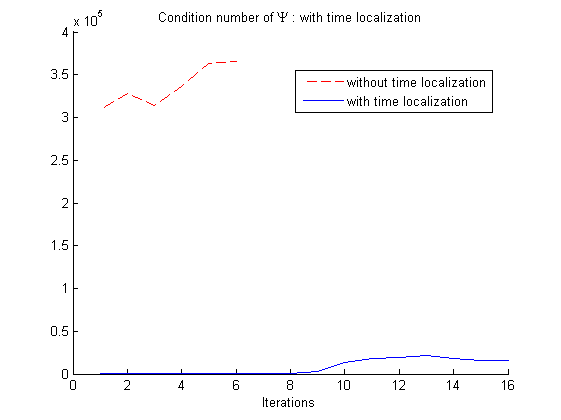}
\includegraphics*[width=5cm]{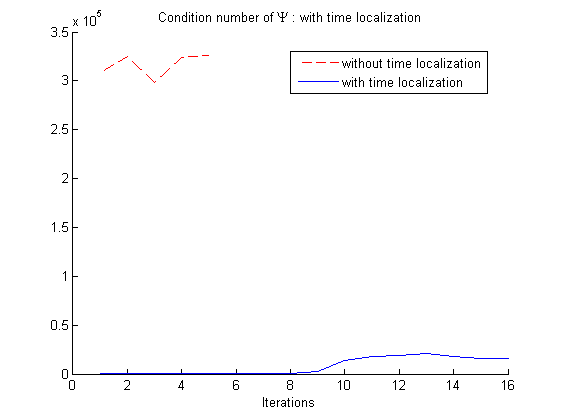}
\includegraphics*[width=5cm]{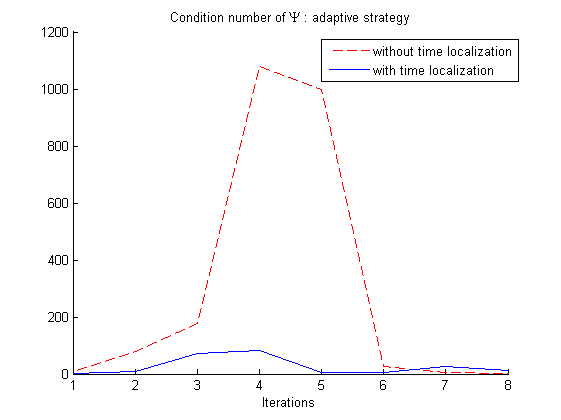}
\includegraphics*[width=5cm]{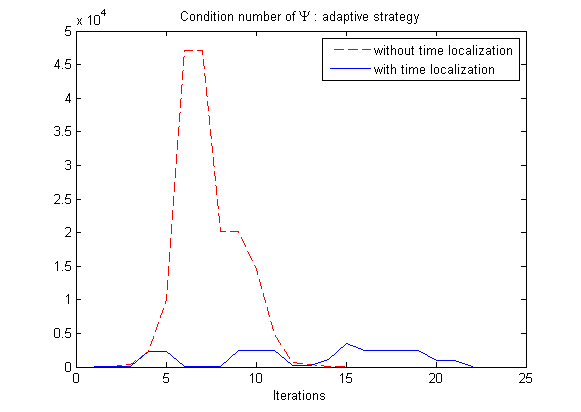}
\end{center}
\caption{\small \textit{Condition number of the sensitivity matrix with (blue line) and without (red dotted line) time localization. First row: finest subdivision. Second row: adaptive parametrization. Left: test 3. Right: test 9.}}\normalsize
\label{cond_adattativa}
\end{figure}

Moreover, as analyzed in section \ref{fiume_1d}, at the stationary regime, the problem becomes ill-conditioned: thus, considering only the transitional regime, time localization could limit the ill-conditioning of the problem. This is evident e.g. in figure \ref{cond_adattativa}, where the four algorithms are compared: without time localization (red dotted line) the condition number of the sensitivity matrix has a much higher upper bound. 

\subsection{Sensitivity of the fourth algorithm to thresholds variations}
It is interesting to analyze what happens when thresholds used in the fourth algorithm are changed. $\epsilon_1$ decides when a the segment corresponding to a parameter should be refined: it is important to keep it not too low, to avoid over-refinements. $\epsilon_2$ is such that parameters less than it are not considered to build the sensitivity matrix: avoiding small parameters reduces computational cost and the ill-conditioning of the problem, since we expect that they are not effective in output variations. 

Previous observations are summarized in table \ref{ese1_soglie}, where test 1 is considered to understand how convergence results varies when thresholds are slightly changed: when $\epsilon_1$ is decreased the over-refinement increases, while when $\epsilon_2$ is lower both the computational cost (number of iterations) and the condition number increase. When both $\epsilon_1$ and $\epsilon_2$ decrease both the distance from the optimal subdivision and the computational cost and the condition number increase. Thus in general to reduce the cost is it better to increase $\epsilon_1$, while to obtain more accurate results it could be useful to adopt smaller $\epsilon_1$ and $\epsilon_2$. 
\begin{center}
\begin{table}[t]
\tiny{\begin{tabular}{|c|c||c|c|c|c|c|c|c|} \hline
\multicolumn{1}{|c|}{$\epsilon_1$}&\multicolumn{1}{|c||}{$\epsilon_2$}&\multicolumn{2}{|c|}{$L^1$ error:} & \multicolumn{2}{|c|}{opt. sub.:}& \multicolumn{1}{|c|}{$\mathcal{F}_d(\mbox{\boldmath $\vartheta$})$} & \multicolumn{1}{|c|}{num. it.}& \multicolumn{1}{|c|}{mean condition number of $\Psi$}\\ 
\multicolumn{1}{|c|}{} &\multicolumn{1}{|c||}{} & \multicolumn{1}{|c|}{up} & \multicolumn{1}{|c|}{down} & \multicolumn{1}{|c|}{up} & \multicolumn{1}{|c|}{down} & \multicolumn{1}{|c|}{} & \multicolumn{1}{|c|}{} & \multicolumn{1}{|c|}{}\\ \hline
0.4 & 0.4 &  1.15 & 0.168 & +1 & 0 & $10^{-4}$ & 7 & 79.9513\\ \hline\hline
0.3 & 0.4 &  1.15 & 0.168 & +1 & +1 & $10^{-4}$ & 7 & 79.9513\\ \hline
0.01 & 0.4 &  1.15 & 0.168 & +1 & +7  & $10^{-5}$ & 7 & 79.9513\\ \hline\hline
0.4 & 0.3 &  1.192 & 0.02 & +1 & 0 & $10^{-5}$ & 8 & 173.2498\\ \hline
0.4 & 0.01 &  1.207 & 0.05 & +1 & 0 & $10^{-6}$ & 9 & 252.7891\\ \hline\hline
0.01 & 0.01 &  1.259 & 0.01 & +1 & +3 & $10^{-6}$ & 9 & 210.4405\\ \hline
\end{tabular}}
\caption{\small \textit{Test 1: results for different values of $\epsilon_1$ and $\epsilon_2$. 
}}  
\label{ese1_soglie}
\end{table}
\end{center}

\section{Conclusions}
This paper presents a mathematical algorithm to solve a class of parabolic inverse problems based upon a convection-diffusion-reaction equation, extending some ideas presented in \cite{deolmi} and \cite{marcuzzi2}. Both liquid (e.g. water) and gas (e.g. air) pollution problems could be considered: when source location is known, we have demonstrated that the problem is well-posed and can be solved e.g. using the Projected damped Gauss Newton method. When $\Gamma_{in}$ is unknown, we have proved that adaptive parametrization with time localization is an effective strategy to estimate a sparse vector of parameters. 

\begin{figure}[h]
\begin{center}
\includegraphics*[width=7cm]{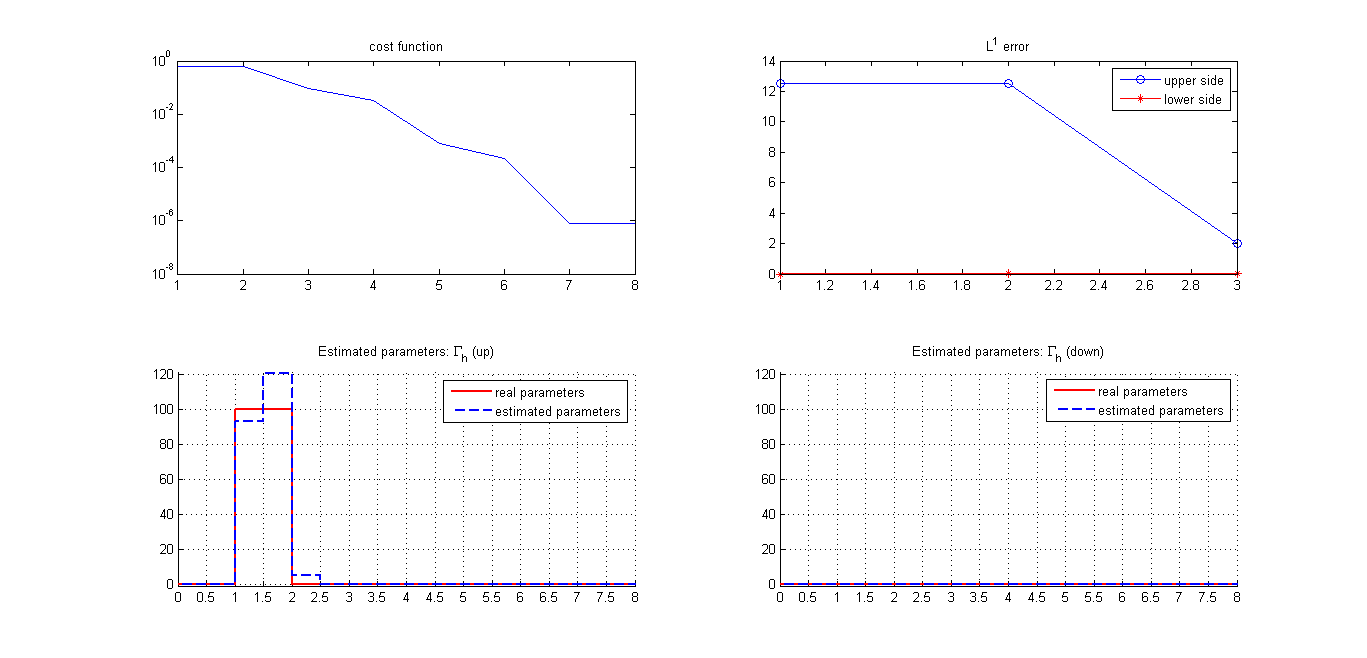}
\end{center}
\caption{\small \textit{Need of an under-refinement strategy.}}\normalsize
\label{2param_unref}
\end{figure}

It could be interesting to introduce also an unrefinement strategy, trying to get closer to the optimal subdivision. For example consider figure \ref{2param_unref}: the optimal strategy would estimate only one parameter in $[1,2]\times \left\{1\right\}$, and it would not bisect the segment $[1,2]$. Instead algorithm 4 bisects $[1,2]$: the problem here is that the direction of the convective field $\textbf{u}$ produces an overestimate of the right hand side parameter of $[1,2]$ and an underestimate of the left hand side one.

Another interesting aspect could be the generalizzation of the problem to time varying boundary conditions on $\Gamma_{in}$ and to analyze more deeply the problem when space-time varying velocity fields are considered.

\section*{Acknowledgments}
The authors are grateful to P. Mannucci for helpful discussions and her fundamental contribution in the demonstration of Lemma \ref{monotonia}.
\appendix
\section{The importance of stabilizing the problem}
Dealing with convection dominated problems ($\left\|\textbf{u}\right\|>>\mu$) could be problematic, due to spurious oscillations caused by the standard FE method. The simplest way to stabilize the problem is to refine the mesh, i.e. to consider a higher number of degrees of freedom; otherwise on a coarse mesh a stabilization method such as SUPG, DW or GLS, to mention only some of them, should be used. To simplify the problem in the following we apply the simplest strategy, i.e. we refine the mesh. However a stabilization method could be included in the model, modifying the weak FE formulation. Stabilization techniques are used e.g. in \cite{becher,collis}.

\begin{figure}[t]
\begin{center}
\includegraphics*[width=5cm]{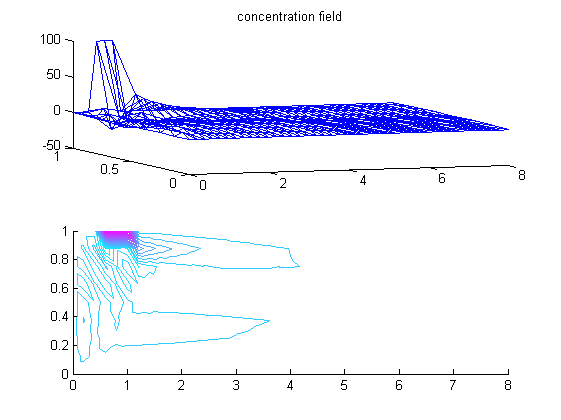}
\includegraphics*[width=7cm]{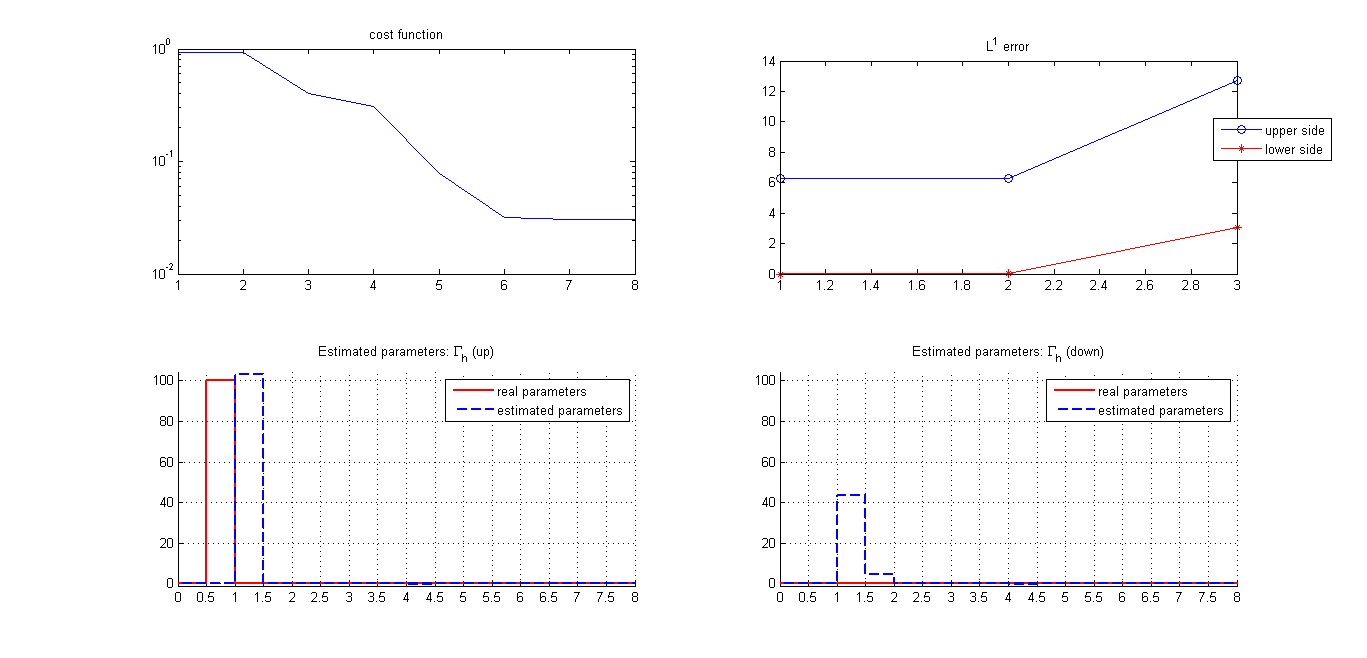}
\includegraphics*[width=5cm]{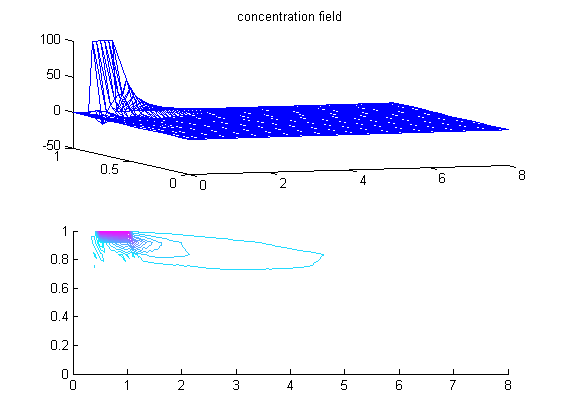}
\includegraphics*[width=7cm]{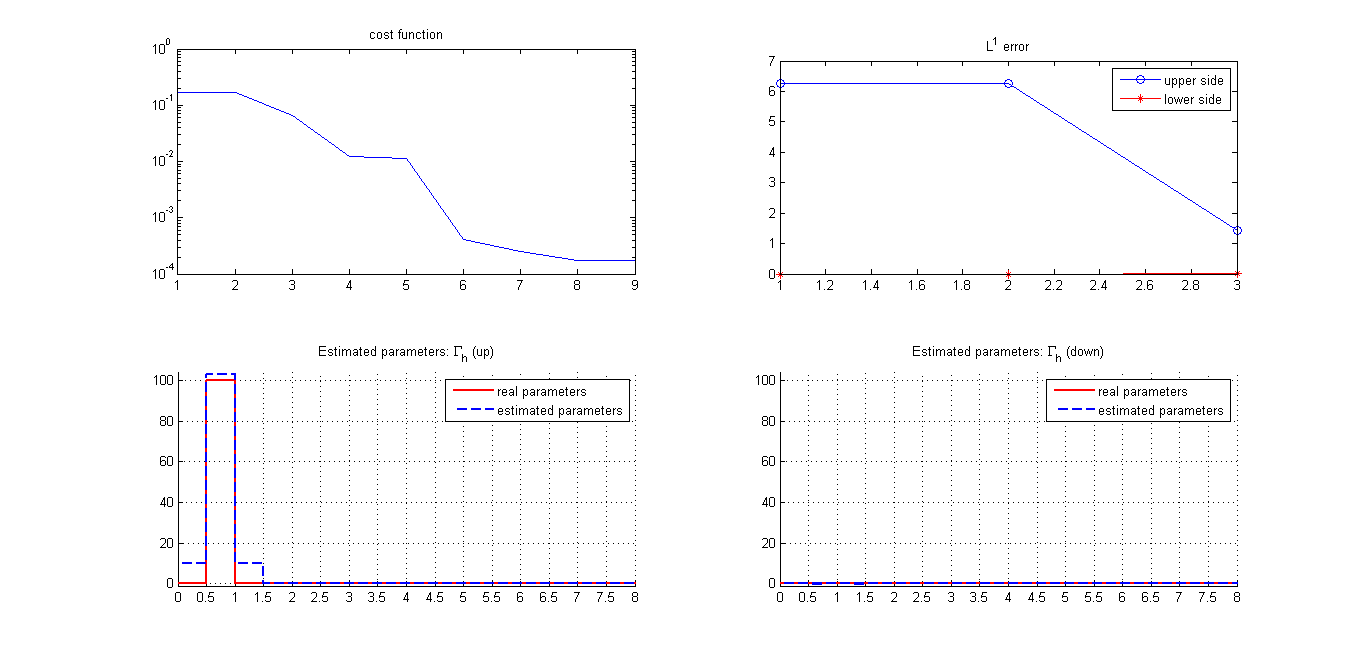}
\includegraphics*[width=5cm]{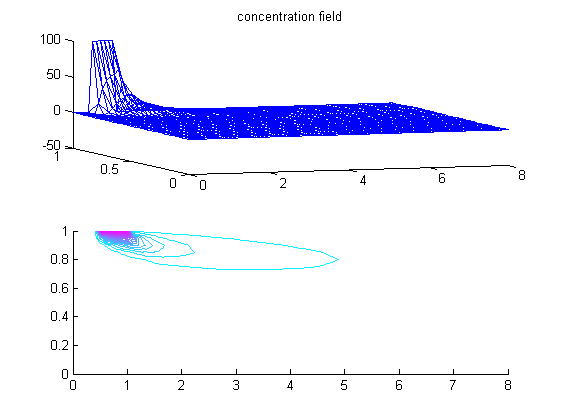}
\includegraphics*[width=7cm]{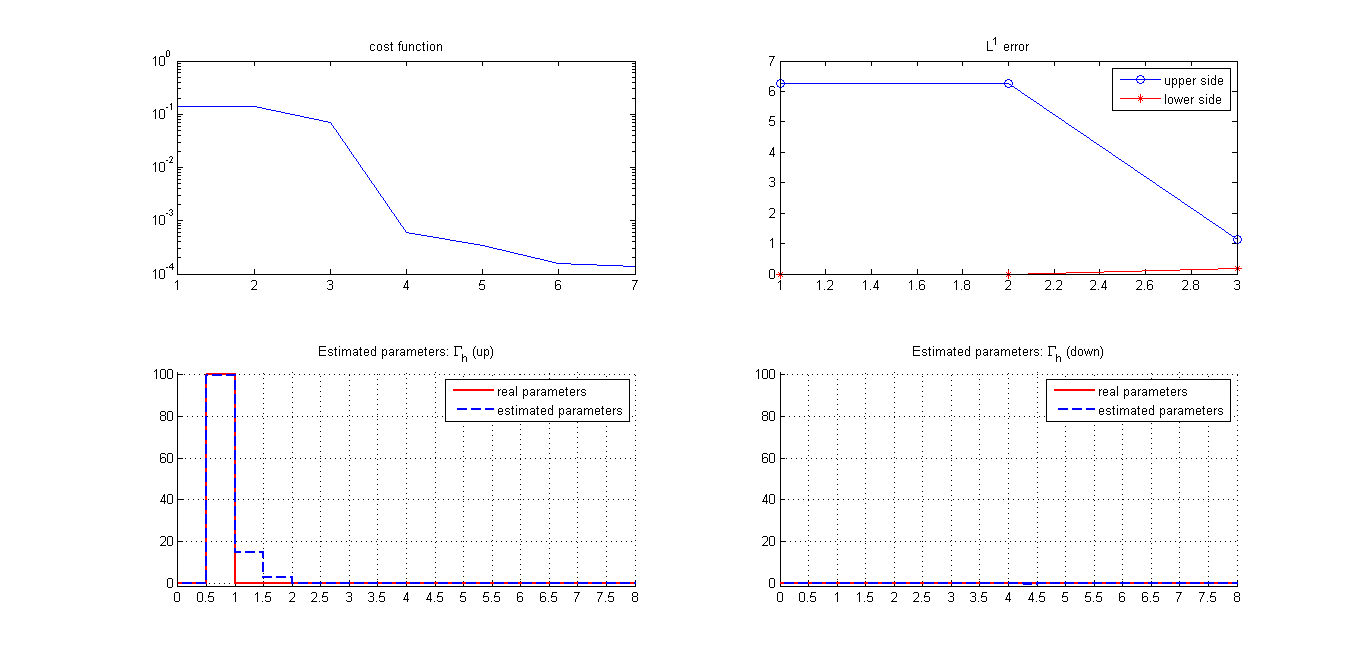}
\end{center}
\caption{\small \textit{Importance of using stabilization: concentration field (left), estimated profile (right). First row: using $41$ nodes along x-axis and $9$ along y-axis. Second row: using $81$ nodes along x-axis and $13$ along y-axis. Third row: using $81$ nodes along x-axis and $21$ along y-axis.}}\normalsize
\label{stab_fiume}
\end{figure}

In this section we want to point out that the problem must be stabilized to obtain a correct estimate. In fact consider $\Omega = [0,8]\times [0,1]$,  $\Gamma_h=[0,8]\times\left\{1\right\}\cup [0,8]\times\left\{0\right\}$, the velocity field $\textbf{u}$ is modelled as a Poiseuille flow i.e.
$$\textbf{u}(x_1,x_2)=\left(\begin{array}{c c} -4\nu x^2_2 + 4 \nu x_2\\ 0 \end{array}\right),$$
assume moreover that $\mu=0.1$, $\sigma=0.1$, $C_{up}=0.1$ and $\Gamma_{in}=[0.5,1]\times\left\{1\right\}$, $\vartheta = 100$. Apply to it the adaptive strategy with time localization, on different meshes. Results are depicted in figure \ref{stab_fiume}.
As it can be seen, when the mesh is too coarse, the presence of spurious oscillations compromise the convergence of the algorithm to the real profile, whereas adopting a fine mesh eliminates them and gives a good estimate of the boundary control.

\end{document}